\documentclass[a4paper,12pt]{article}

\usepackage{amsmath, amsthm, amsfonts, amssymb}
\usepackage{graphicx}

\theoremstyle{plain}
\newtheorem{theorem}{Theorem}[section]

\newtheorem{lemma}{Lemma}[section]

\newtheorem{corollary}[theorem]{Corollary}

\theoremstyle{remark}
\newtheorem{remark}{Remark}[section]

\numberwithin{equation}{section}

\def\tht{\theta}

\def\e{\varepsilon}
\def\g{\gamma}

\def\l{\lambda}
\def\p{\partial}

\def\k{\varkappa}
\def\E{\mbox{\rm e}}
\def\a{\alpha}
\def\b{\beta}

\def\d{\delta}

\def\z{\zeta}

\def\vp{\varphi}

\def\Odr{\mathcal{O}}
\def\H{W_2}
\def\Hloc{W_{2,loc}}

\def\di{\,\mathrm{d}}

\def\I{\mathrm{I}}
\def\iu{\mathrm{i}}

 \DeclareMathOperator{\RE}{Re}
\DeclareMathOperator{\IM}{Im} \DeclareMathOperator{\spec}{\sigma}

\DeclareMathOperator{\conspec}{\sigma_{c}}
\DeclareMathOperator{\pointspec}{\sigma_{p}}
\DeclareMathOperator{\resspec}{\sigma_{r}}
\DeclareMathOperator{\dist}{dist}
\DeclareMathOperator{\supp}{supp} \DeclareMathOperator{\sgn}{sgn}

\begin{document}

\allowdisplaybreaks


\title{\textbf{On spectrum of a periodic operator with a small
localized perturbation}}
\author{D.~Borisov$^{a,b}$ and R.~Gadyl'shin$^{b}$}
\date{}
\maketitle

\footnote{The work is supported in parts by RFBR (06-01-00138 --
D.B., 05-01-97912-r\_agidel -- R.G.). D.B. is also supported by
\emph{Marie Curie International Fellowship} within 6th European
Community Framework (MIF1-CT-2005-006254).}

\vspace{-1 true cm}

\begin{quote}
{\small {\em a) Nuclear Physics Institute, Academy of Sciences,
25068 \v Re\v z
\\
\phantom{a) }near Prague, Czechia
\\
b) Bashkir State Pedagogical University, October Revolution
St.~3a,
\\
\phantom{a) }450000 Ufa, Russia
 \\
E-mails: \texttt{borisovdi@yandex.ru},
\texttt{gadylshin@yandex.ru}}}
\end{quote}

\begin{abstract}
We study the spectrum of a periodic self-adjoint operator on the
axis perturbed by a small localized nonself-adjoint operator. It
is shown that the continuous spectrum is independent of the
perturbation, the residual spectrum is empty, and the point
spectrum has no finite accumulation points. We address the
existence of the embedded eigenvalues. We establish the
necessary and sufficient conditions of the existence of the
eigenvalues and construct their asymptotics expansions. The
asymptotics expansions for the associated eigenfunctions are
also obtained. The examples are given.
\end{abstract}

\section{Introduction}

It is well-known that the spectrum of a self-adjoint periodic
one-dimensional differential operator consists of zones
separated by lacunas (see, for instance, (\cite[Ch. V, Sec.
56]{Gl}, \cite[Ch. 5]{E}). Perturbation of such operator by a
rapidly decreasing potential does not change the continuous part
of the spectrum but produces the isolated eigenvalues in the
lacunas. The existence and the number of such eigenvalues were
studied, for instance, in  \cite{RB1}, \cite{Zl1}, \cite{Zl2},
\cite{GR}. It was shown that the number of the eigenvalues in
each lacuna is finite and there are at most two eigenvalues in
the distant lacunas. In \cite{GR} they also studied the case
when the perturbing potential is multiplied by a small coupling
constant. It was established that each lacuna contains at most
two eigenvalues. The necessary and sufficient conditions exactly
determining the number of the eigenvalues in a given lacuna were
adduced.

In the present paper we study the spectrum of a self-adjoint
periodic second-order differential operator on the axis
perturbed by a linear operator of the form $\e \mathcal{L}_\e$,
where $\e$ is a small positive parameter. The main feature of
$\mathcal{L}_\e$  is that it is localized in the following
sense. The support of the function $\mathcal{L}_\e u$ lies
inside some fixed finite segment  $\overline{Q}$ and this
function is fully determined by the values the argument $u$
takes on $Q$.

The main difference of the perturbation we study from the cases
in the papers cited is that we do not assume the
self-adjointness neither for $\mathcal{L}_\e$ nor for the
perturbed operator. Moreover, the set of possible perturbation
described by such perturbation includes, apart from the
potentials, a wide class of examples of various nature like
differential operator, integral operator, linear functional (see
last section).

In the paper we show that the continuous spectrum of the
perturbed operator coincides with the spectrum of the
unperturbed one. We also establish that the residual spectrum is
empty, while the point spectrum consists of at most countably
many eigenvalues of finite multiplicity and has no finite
accumulation points. We also give an example of the perturbation
which originates an embedded eigenvalue. We notice that similar
phenomenon could not rise in the problems considered in
\cite{RB1}, \cite{Zl1}, \cite{Zl2}, \cite{GR}. We also provide
the sufficient conditions guaranteeing absence of the embedded
eigenvalues. It is established that the perturbed eigenvalues
tend either to infinity or to the edges of non-degenerated
lacunas in the continuous spectrum. We prove that there exists
at most one such eigenvalue in the vicinity of an edge of a
given non-degenerated lacuna. We give the criteria for the
existence of this eigenvalue and construct its asymptotics
expansions as well as the expansion for the associated
eigenfunction.

In conclusion let us describe briefly the structure of the
paper. In the following section we formulate the problem and
present the main results. In the third section we prove the
general theorem on the position of the perturbed spectrum and
show that continuous spectrum is independent of the perturbation
and the residual one is empty. In the fourth section we show the
absence of the embedded eigenvalues in the finite part of
spectrum if $\e$ is small enough. The fifth section is devoted
to the countability, convergence and some other properties of
the point spectrum. In the sixth section we study the existence
of the embedded eigenvalues. In the seventh section certain
auxiliary statements are proven. These statements are employed
in the eighth section where we construct the asymptotics for the
eigenvalues converging to the edges of the non-degenerate
lacunas in the continuous spectrum. In the last ninth section we
give examples of the perturbation and apply to them the general
results of the work.

\section{Formulation of the problem and the main results}

Let
\begin{equation*}
\mathcal{H}_0:=-\frac{d}{dx}p\frac{d}{dx}+q
\end{equation*}
be a self-adjoint operator in $L_2(\mathbb{R})$ with the domain
$\H^2(\mathbb{R})$. Here  $p=p(x)$ is $1$-periodic piecewise
continuously differentiable real function,  $q=q(x)$ is
$1$-periodic piecewise continuous real function, and
\begin{equation}\label{0.1}
p(x)\geqslant p_0>0, \quad x\in \mathbb{R}.
\end{equation}
Without loss of generality throughout the paper we assume that
$p(0)=1$.

Let $\mathcal{L}_\e: \H^2(Q)\to L_2(Q)$ be a linear operator
bounded uniformly in $\e$, and generally speaking unbounded as
an operator in $L_2(Q)$. We introduce the operator mapping
$\Hloc^2(\mathbb{R})$ into $L_2(\mathbb{R})$ by the following
rule: an element from $\Hloc^2(\mathbb{R})$ is restricted to
$Q$, then the operator $\mathcal{L}_\e$ is applied, and the
result is continued by zero outside $Q$. Such operator is
naturally to indicate by $\mathcal{L}_\e$. Clearly, this is an
unbounded operator in $L_2(\mathbb{R})$ with the domain
$\Hloc^2(\mathbb{R})$.

We indicate $\mathcal{H}_\e:= (\mathcal{H}_0-\e \mathcal{L}_\e)$
considering it as an operator in $L_2(\mathbb{R})$ having
$\H^2(\mathbb{R})$ as the domain. The operator $\mathcal{H}_\e$
is closed (see Lemma~\ref{lm2.1}).

The main aim of this paper is to study the behaviour of the
spectrum of the operator $\mathcal{H}_\e$ as $\e\to0$. Before
presenting the main results we introduce additional notations
and remind some known facts.

We will employ the symbols $\spec(\cdot)$, $\conspec(\cdot)$ and
$\pointspec(\cdot)$ to indicate the spectrum, continuous
spectrum, and the point spectrum, while
\begin{equation*}
\resspec(\cdot):=\spec(\cdot)\backslash
\left(\conspec(\cdot)\cup\pointspec(\cdot)\right)
\end{equation*}
is the residual spectrum. It is known \cite[Ch. 2, Sec.
2.2,\,2.3, Ch. 5, Sec. 5.3]{E} that the operator $\mathcal{H}_0$
has a band spectrum
\begin{equation}\label{2.1a}
\spec(\mathcal{H}_0)=\conspec(\mathcal{H}_0)=\bigcup_{n=0}^\infty
\big[\mu_n^+,\mu_{n+1}^-\big],
\end{equation}
where the numbers
\begin{equation*}
\mu_0^+<\mu_1^-\leqslant\mu_1^+<\mu_2^-\leqslant\mu_2^+
<\mu_3^-\leqslant\mu_3^+<\ldots
\end{equation*}
are simple eigenvalues of the boundary value problems
\begin{equation}
\begin{gathered}
\left(-\frac{d}{dx}p\frac{d}{dx}+q\right) \phi_n^\pm=
\mu_n^\pm\phi_n^\pm,\quad x\in(0,1),
\\
\phi_n^\pm(0)+(-1)^{n+1}\phi_n^\pm(1)=0,\quad
\frac{d\phi_n^\pm}{dx}(0)+(-1)^{n+1}\frac{d\phi_n^\pm}{dx}(1)=0.
\end{gathered}\label{1.3}
\end{equation}
For $a\in\mathbb{C}$, $\d>0$ we denote $S_\d(a):=\{\l\in
\mathbb{C}: |\arg (\l-a)|<\d\}$.

We are ready to formulate the main results of the paper.

\begin{theorem}\label{th2.2}
There exist positive $\d_i=\d_i(\e)\xrightarrow[\e\to0]{}0$,
$i=1,2$, so that for $\e$ small enough the inclusion
$\spec(\mathcal{H}_\e)\subset S_{\d_2(\e)}\big(\mu_0^+ -
\d_1(\e)\big)$ holds true.
\end{theorem}

\begin{theorem}\label{th2.1}
For $\e$ small enough the identities
$\conspec(\mathcal{H}_\e)=\conspec(\mathcal{H}_0)$,
$\resspec(\mathcal{H}_\e)=\emptyset$ are valid.
\end{theorem}

\begin{theorem}\label{th2.33}
The point spectrum of the operator $\mathcal{H}_\e$ consists of
countably many eigenvalues of finite multiplicities and has no
finite accumulation points.
\end{theorem}

\begin{theorem}\label{th2.3}
Let $K$ be a compact set in the complex plane such that
$K\cap\conspec(\mathcal{H}_0)\not=\emptyset$. Then for $\e$
small enough the set $\conspec(\mathcal{H}_\e)\cap K$ contains
no embedded eigenvalues.
\end{theorem}

We stress that this theorem does not exclude the presence of the
embedded eigenvalues tending to infinity as $\e\to0$. In the
sixth section we will give an example of the operator
$\mathcal{H}_\e$ which has an embedded eigenvalue. In the
following theorem we provide the sufficient conditions of the
absence of such eigenvalues.

\begin{theorem}\label{th2.31}
Assume that at least one of the following conditions is valid
\begin{enumerate}\def\theenumi{(\arabic{enumi})}

\item\label{it1th2.4}
For any subinterval $\widetilde{Q}\subseteq Q$ the estimate
\begin{equation}\label{1.3e}
\|\mathcal{L}_\e u\|_{L_2(\widetilde{Q})}\leqslant
C\|u\|_{\H^2(\widetilde{Q})}
\end{equation}
holds true with the constant $C$ independent of $\e$ and
$\widetilde{Q}$.

\item\label{it2th2.4} The operator $\mathcal{L}_\e$ can be represented
as
\begin{equation*}
\mathcal{L}_\e=\frac{d}{dx}a_\e\frac{d}{dx}+
\widetilde{\mathcal{L}}_\e,
\end{equation*}
where $a_\e$ is piecewise continuously differentiable function
having support inside $\overline{Q}$ and satisfying the relation
\begin{equation}\label{1.3d}
\e\max\limits_{\overline{Q}} |a_\e'(x)|\xrightarrow[\e\to0]{}0,
\end{equation}
and $\widetilde{\mathcal{L}}_\e:\H^1(Q)\to L_2(Q)$ is a linear
operator bounded uniformly in $\e$.
\end{enumerate}
Then for $\e$ small enough the continuous spectrum of
$\mathcal{H}_\e$ contains no embedded eigenvalues.
\end{theorem}

\begin{theorem}\label{th2.32}
Let $K$ be a compact set in the complex plane such that
$K\cap\conspec(\mathcal{H}_0)\not=\emptyset$. Then for $\e$
small enough each of the eigenvalues $\mathcal{H}_\e$ not
leaving $K$ for all $\e$ small enough converges in the limit
$\e\to0$ to one of the edges of a non-degenerate lacuna in the
part of the spectrum of $\mathcal{H}_0$ lying inside $K$.
\end{theorem}

Let $\tht_i(x,\lambda)$ be the solutions to the equation
\begin{align}\label{1.4}
&\left(-\frac{d}{dx}p\frac{d}{dx}+q-\l\right)v=0,\quad x\in
\mathbb{R} ,
\end{align}
satisfying the initial conditions
\begin{equation}\label{1.3a}
\tht_1(0,\lambda)=1,\quad \frac{d\tht_1}{dx}(0,\lambda)=0,
\qquad\tht_2(0,\lambda)=0,\quad \frac{d\tht_2}{dx}(0,\lambda)=1,
\end{equation}
where $\l$ is a complex parameter. For the sake of brevity
hereinafter we denote $\tht_i(\l):=\tht_i(1,\l)$,
$\tht'_i(\l):=\frac{d\tht_i}{dx} (1,\l)$, $i=1,2$. We set
$D(\l):=\tht_1(\l)+\tht'_2(\l)$.

Let $\mu_n^\pm$ be one of the edges of a non-degenerate lacuna
in the spectrum of $\mathcal{H}_0$. We choose the eigenfunctions
of the problem (\ref{1.3}) being real and continue them
1-periodically over the axis for even $n$ and 1-antiperiodically
for odd $n$. It is clear that the functions continued are twice
piecewise continuously differentiable. We normalize them as
follows
\begin{equation}\label{1.3c}
|\phi_{n}^\pm(0)|^2+\left|\frac{d\phi_{n}^\pm}{dx}(0)\right|^2
=|\tht_1'(\mu_n^\pm)|+|\tht_2(\mu_n^\pm)|.
\end{equation}
We will show below that the right hand side of this identity is
non-zero (see Item~\ref{it4lm1.1} of Lemma~\ref{lm7.3}), and
this is why this normalization makes sense.

Let $\mathcal{G}_{n,0}^\pm$ be an integral operator defined on
$L_2(Q)$:
\begin{gather}
(\mathcal{G}_{n,0}^\pm f)(x):=\int\limits_\mathbb{R}
G_{n,0}^\pm(x,t)f(t)\di t,\label{2.8a}
\\
G_{n,0}^\pm(x,t):=\frac{1}{2}\left\{
\begin{aligned}
&\tht_1(t,\mu_n^\pm)\tht_2(x,\mu_n^\pm)-
\tht_1(x,\mu_n^\pm)\tht_2(t,\mu_n^\pm),&&t>x,
\\
&\tht_1(x,\mu_n^\pm)\tht_2(t,\mu_n^\pm)-
\tht_1(t,\mu_n^\pm)\tht_2(x,\mu_n^\pm),&&t<x.
\end{aligned} \right.\nonumber
\end{gather}
Since $\mathcal{G}_{n,0}^\pm:L_2(Q)\to\H^2(Q)$ is a bounded
operator, it follows that for $\e$ small enough the operator
$\mathcal{L}_\e\mathcal{G}_{n,0}^\pm: L_2(Q)\to L_2(Q)$ is
bounded uniformly in $\e$, and thus for $\e$ small enough the
bounded operator
\begin{equation*}
\mathcal{A}_n^\pm(\e,0):=\left(\I-\e
\mathcal{L}_\e\mathcal{G}_{n,0}^\pm\right)^{-1}.
\end{equation*}
is well-defined in $L_2(Q)$. Hereinafter  $\I$ is the identity
mapping. By the dot we will indicate the differentiation w.r.t.
$\lambda$.

\begin{theorem}\label{th2.5}
Let $\mu_n^\pm$ be one of the edges of a non-degenerate lacuna
in the spectrum of  $\mathcal{H}_0$. Then the operator ${H}_\e$
has at most one eigenvalue $\l_{\e,n}^\pm$ converging to
$\mu_n^\pm$ as $\e\to0$. This eigenvalue exists if and only if
\begin{equation}\label{2.9}
\pm\RE \big(\phi_{n}^\pm,
\mathcal{A}_n^\pm(\e,0)\mathcal{L}_\e\phi_{n}^\pm\big)_{L_2(Q)}>0.
\end{equation}
If exists, this eigenvalue is simple and has the asymptotics
expansion
\begin{align}
&\l_{\e,n}^\pm=\mu_n^\pm\mp\frac{\e^2}{4|\overset{\,\textbf{.}}{D}(\mu_n^\pm)|}
\big(\phi_{n}^\pm,
\mathcal{A}_n^\pm(\e,0)\mathcal{L}_\e\phi_{n}^\pm\big)_{L_2(Q)}^2
(1+\Odr(\e^2)),\label{2.10}
\\
&\l_{\e,n}^\pm=\mu_n^\pm\mp \e^2 \big(k_{n,\e}^{\pm,1}+\e
k_{n,\e}^{\pm,2}\big)^2+\Odr(\e^4
|k_{n,\e}^{\pm,1}|+\e^5),\label{2.11}
\\
&k_{n,\e}^{\pm,1}:=\pm\frac{\big(\mathcal{L}_\e\phi_{n}^\pm,
\phi_{n}^\pm\big)_{L_2(Q)}}{2\sqrt{|\overset{\,\textbf{.}}{D}(\mu_n^\pm)|}},\quad
k_{n,\e}^{\pm,2}:=\pm\frac{\big(\mathcal{L}_\e
\mathcal{G}_{n,0}^\pm\mathcal{L}_\e \phi_{n}^\pm,\phi_{n}^\pm
\big)_{L_2(Q)}}{2\sqrt{|\overset{\,\textbf{.}}{D}(\mu_n^\pm)|}}.
\label{2.8b}
\end{align}
The asymptotics expansion for the associated eigenfunction reads
as follows
\begin{equation}\label{2.13}
\psi_{\e,n}^\pm=\phi_{n}^\pm+\e \mathcal{G}_{n,0}^\pm
\mathcal{L}_\e \phi_{n}^\pm +\Odr(\e^2)
\end{equation}
in the norm of $\H^2(\a_1,\a_2)$ for any $\a_1,\a_2\in
\mathbb{R}$.
\end{theorem}

\begin{remark}\label{rm1.2}
In Lemma~\ref{lm7.3} we will give the formulas for
$\overset{\,\textbf{.}}{D}(\mu_n^\pm)$, which imply in
particular that $\overset{\,\textbf{.}}{D}(\mu_n^\pm)\not=0$, if
$\mu_n^\pm$ is an edge of a non-degenerate lacuna.
\end{remark}

Theorems~\ref{th2.32},~\ref{th2.5} yield immediately

\begin{corollary}\label{cr1.1}
Let $K$ be a compact set in the complex plane. Then each of the
eigenvalues of $\mathcal{H}_\e$ lying inside $K$ for all $\e$
small enough is simple.
\end{corollary}

\begin{theorem}\label{th2.6}
Let $\mu_n^\pm$ be one of the edges of a non-degenerate lacuna
in the spectrum of  $\mathcal{H}_0$. If
\begin{equation}\label{1.12a}
\RE \left(k_{n,\e}^{\pm,1}+\e k_{n,\e}^{\pm,2}\right)\geqslant
C(\e)\e^2, \quad C(\e)\xrightarrow[\e\to0]{}+\infty,
\end{equation}
there exists the eigenvalue $\l_{\e,n}^\pm$, and the identities
(\ref{2.10}), (\ref{2.11}) hold true. In the case
\begin{equation}\label{1.12b}
\RE \left(k_{n,\e}^{\pm,1}+\e k_{n,\e}^{\pm,2}\right)\leqslant
-C(\e)\e^2, \quad C(\e)\xrightarrow[\e\to0]{}+\infty,
\end{equation}
the operator $\mathcal{H}_\e$ has no eigenvalues converging to
$\mu_n^\pm$ as $\e\to0$.
\end{theorem}

\begin{remark}\label{rm1.1}
We also obtain the explicit formula for the eigenfunction
$\psi_{\e,n}^\pm$ and describe how it behaves at infinity (see
(\ref{2.12}), (\ref{2.14}), (\ref{2.15})).
\end{remark}

\begin{remark}\label{rm1.3}
We notice that in the particular case $p\equiv const$, $q\equiv
const$ the continuous spectrum of $\mathcal{H}_0$ coincides with
the semi-axis $[q,+\infty)$ and has no internal lacunas. In this
case the semi-infinite lacuna $(-\infty,q)$ is the only
non-degenerate one and
Theorems~\ref{th2.32},~\ref{th2.5},~\ref{th2.6} describes the
behaviour of the eigenvalues in the vicinity of the point
$\mu_0^+=q$.
\end{remark}

\section{
Proof of Theorems~\ref{th2.2},~\ref{th2.1} }

We denote $B_r(a):=\{\l\in \mathbb{C}: |\l-a|<r\}$. For any pair
of non-empty sets $M_1,M_2\subset\mathbb{C}$ we indicate
\begin{align*}
\dist(M_1,M_2):=\inf\limits_{\genfrac{}{}{0 pt} {}{\l_1\in
M_1}{\l_2\in M_2}}|\l_1-\l_2|.
\end{align*}

\begin{lemma}\label{lm3.1}
Let $M$ be a non-empty closed set in the complex plane such that
$M\cap\spec(\mathcal{H}_0)=\emptyset$ and for some $a\in
\mathbb{R}$, $\d\in(0,\pi/2)$, $r>0$ the inclusion $M\backslash
B_r(0) \subset \mathbb{C}\setminus S_\d(a)$ is valid. Then for
all $f\in L_2(\mathbb{R})$ and $\l\in M$ the estimate
\begin{equation*}
\|(\mathcal{H}_0-\l)^{-1}f\|_{\H^2(\mathbb{R})}\leqslant C
\|f\|_{L_2(\mathbb{R})},
\end{equation*}
is valid where the constant $C$ is independent of $\l\in M$.
\end{lemma}

\begin{proof}
Let $f\in L_2(\mathbb{R})$, $\l\in M$. Since
$\l\not\in\spec(\mathcal{H}_0)$, it follows that the operator
$(\mathcal{H}_0-\l)^{-1}: L_2(\mathbb{R})\to L_2(\mathbb{R})$ is
bounded and in accordance with formula (3.16) in \cite[Ch. V,
Sec. 3.5]{K} the estimate
\begin{equation*}
\|u\|_{L_2(\mathbb{R})}\leqslant
\frac{\|f\|_{L_2(\mathbb{R})}}{\dist(\lambda,\spec(\mathcal{H}_0))},
\end{equation*}
holds true for all $\l\in M$, where
$u:=(\mathcal{H}_0-\l)^{-1}f$. The inequality (\ref{0.1}) and
the boundedness of the function $q$ imply
\begin{align*}
&(pu',u')_{L_2(\mathbb{R})}+
(qu,u)_{L_2(\mathbb{R})}-\l\|u\|^2_{L_2(\mathbb{R})}=
(f,u)_{L_2(\mathbb{R})},
\\
&\|u'\|^2_{L_2(\mathbb{R})}\leqslant
C(1+|\l|)\|u\|^2_{L_2(\mathbb{R})}+\|f\|_{L_2(\mathbb{R})}
\|u\|_{L_2(\mathbb{R})},
\end{align*}
where the constant $C$ is independent of $M$. Now we express the
second derivative of $u$ by the equation $(\mathcal{H}_0-\l)u=f$
and in view of two last inequalities obtain the estimates
\begin{equation}\label{2.23}
\|(\mathcal{H}_0-\l)^{-1}f\|_{\H^2(\mathbb{R})}\leqslant C
\frac{1+|\l|}{\dist(\l,\spec(\mathcal{H}_0))}
\|f\|_{L_2(\mathbb{R})}
\end{equation}
with the constant $C$ independent of $M$. If the set $M$ is
bounded, the statement of the lemma follows from the obtained
estimate. In the case the set $M$ is unbounded the statement of
the lemma follows from the estimate (\ref{2.23}) and an obvious
inequality
\begin{equation*}
\sup\limits_{M\setminus B_r(0)}
\frac{1+|\l|}{\dist(\l,\spec(\mathcal{H}_0))} \leqslant
\sup\limits_{\mathbb{C}\setminus S_\d(a)}
\frac{1+|\l|}{\dist(\l,\spec(\mathcal{H}_0))} <\infty.
\end{equation*}
\end{proof}

\begin{lemma}\label{lm2.1}
The operator $\mathcal{H}_\e$ is closed for all $\e$ small
enough.
\end{lemma}

\begin{proof}
Since $(\mu_0^+-1)\not\in\spec(\mathcal{H}_0)$ by (\ref{2.1a}),
Lemma~\ref{lm3.1} with $M=\{\mu_0^+-1\}$ and the uniform in $\e$
boundedness of the operator $\mathcal{L}_\e$ yield
\begin{align*}
\|\mathcal{L}_\e v\|_{L_2(\mathbb{R})}=\|\mathcal{L}_\e
v\|_{L_2(Q)}\leqslant & C\|v\|_{\H^2(\mathbb{R})}\leqslant C \|
(\mathcal{H}_0-\mu_0^++1)v\|_{L_2(\mathbb{R})}
\\
\leqslant& C\left( \| \mathcal{H}_0
v\|_{L_2(\mathbb{R})}+|\mu_0^+-1|
\|v\|_{L_2(\mathbb{R})}\right).
\end{align*}
Hence, the operator $\e \mathcal{L}_\e$ is
$\mathcal{H}_0$-bounded and for $\e$ small enough its
$\mathcal{H}_0$-bound is strictly less than one. By \cite[Ch.
I\!V, Sec. 1.1, Thm. 1.1]{K} it completes the proof.
\end{proof}

\begin{proof}[Proof of Theorem~\ref{th2.2}]
We choose a pair of numbers $a>0$, $\d\in(0,\pi/2)$. It is
sufficient to show that for $\e$ small enough the inclusion
$\spec(\mathcal{H}_\e)\subset S_\d(\mu_0^+-a)$ is valid. In
turn, this inclusion is equivalent to the existence of the
resolvent of the operator $\mathcal{H}_\e$ for all $\l\in
\mathbb{C}\setminus S_\d(\mu_0^+-a)$ if $\e$ is small enough.
Let us prove the last fact.

The set $M:=\mathbb{C}\setminus S_\d(\mu_0^+-a)$ satisfies the
hypothesis of Lemma~\ref{lm3.1}, and this is why by this lemma
and the uniform boundedness of $\mathcal{L}_\e$ we conclude that
the operator $\mathcal{L}_\e(\mathcal{H}_0-\l)^{-1}$ is bounded
uniformly in $\e$. Thus, for $\e$ small enough the operator
$(\I-\e\mathcal{L}_\e(\mathcal{H}_0-\l)^{-1})$ is boundedly
invertible for all $\l\in M$. Employing this fact it is easy to
check that the resolvent of $\mathcal{H}_\e$ exists for all
$\l\in M$ and is given by the identity
$(\mathcal{H}_\e-\l)^{-1}=(\mathcal{H}_0-\l)^{-1}
\big(\I-\e\mathcal{L}_\e(\mathcal{H}_0-\l)^{-1}\big)^{-1}$.
\end{proof}

Let points $x_0$, $x_1$ be so that
$\overline{Q}\subset(x_0,x_1)$. Without loss of generality we
suppose that $(x_1-x_0)$ is a natural number. By
$\mathcal{H}_0^{(-)}$, $\mathcal{H}_0^{(0)}$,
$\mathcal{H}_0^{(+)}$ we denote the operator
\begin{equation*}
-\frac{d}{dx}p\frac{d}{dx}+q
\end{equation*}
respectively, in $L_2(-\infty,x_0)$, $L_2(x_0,x_1)$,
$L_2(x_1,+\infty)$. As the domains for this operators we choose
the subset of the functions from $\H^2(-\infty,x_0)$,
$\H^2(x_0,x_1)$, $\H^2(x_1,+\infty)$ vanishing, respectively, at
$x_0$, at $x_0$, $x_1$, at $x_1$. Clearly, the operators
$\mathcal{H}_0^{(-)}$, $\mathcal{H}_0^{(0)}$,
$\mathcal{H}_0^{(+)}$ are self-adjoint.

\begin{lemma}\label{lm3.2}
The operator $\big(\mathcal{H}_0^{(0)}-\iu\big)^{-1}:
L_2(x_0,x_1)\to L_2(x_0,x_1)$ is compact.
\end{lemma}
\begin{proof}
Since the operator $\mathcal{H}_0^{(0)}$ is self-adjoint, it
follows that $\spec(\mathcal{H}_0^{(0)})\subset\mathbb{R}$. Thus
the operator $\big(\mathcal{H}_0^{(0)}-\iu\big)^{-1}$ is bounded
as an operator in $L_2(x_0,x_1)$. By Banach theorem on the
inverse operator \cite[Ch. I\!V, Sec. 5.4, Thm. 3]{KF} we
conclude that the operator $(\mathcal{H}_0^{(0)}-\iu)^{-1}:
L_2(x_0,x_1)\to \H^2(x_0,x_1)$ is bounded. Together with the
compactness of the embedding $\H^2(x_0,x_1)$ in $L_2(x_0,x_1)$
it completes the proof.
\end{proof}

\begin{lemma}\label{lm3.3}
The identity
$\conspec\big(\mathcal{H}_0^{(-)}\oplus\mathcal{H}_0^{(+)}\big)=
\conspec(\mathcal{H}_0)$ is valid.
\end{lemma}
\begin{proof}
Let $\H^2(\mathbb{R},x_0)$ be the subset of the functions from
$\H^2(\mathbb{R})$ vanishing at $x_0$, and
$\widetilde{\mathcal{H}}_0$ be the restriction of
$\mathcal{H}_0$ on $\H^2(\mathbb{R},x_0)$. Then the operator
$\mathcal{H}_0$ is a closed extension of
$\widetilde{\mathcal{H}}_0$. It is obvious that the factor-space
$\H^2(\mathbb{R})/\H^2(\mathbb{R},x_0)$ is one-dimensional and
$\mathcal{H}_0$ is hence a finite-dimensional extension of
$\widetilde{\mathcal{H}}_0$. By  \cite[Ch. I, Sec. 1.2, Thm.
4]{Gl} it follows that
\begin{equation}\label{2.24}
\conspec(\widetilde{\mathcal{H}}_0)=\conspec(\mathcal{H}_0).
\end{equation}
The operators $\widetilde{\mathcal{H}}_0$ and
$\mathcal{H}_0^{(-)}\oplus\mathcal{H}_0^{(+)}$ are unitarily
equivalent; the corresponding unitary operator is defined as
\begin{equation*}
\left(\mathcal{U}\big(u^{(-)},u^{(+)}\big)\right)(x):=
\begin{cases}
u^{(-)}(x),& x<x_0,
\\
u^{(+)}(x),& x>x_1,
\end{cases}
\end{equation*}
where $u^{(\pm)}$ belong to the domains of
$\mathcal{H}_0^{(\pm)}$. Therefore,
$\conspec\big(\mathcal{H}_0^{(-)} \oplus\mathcal{H}_0^{(+)}\big)
=\conspec(\widetilde{\mathcal{H}}_0)$, that together with
(\ref{2.24}) completes the proof.
\end{proof}

\begin{proof}[Proof of Theorem~\ref{th2.1}]
Let $\mathcal{H}_\e^{(0)}:=\mathcal{H}_0^{(0)}-\e
\mathcal{L}_\e$ be an operator in $L_2(x_0,x_1)$ whose domain
coincides with one of $\mathcal{H}_0^{(0)}$. Here the operator
$\mathcal{L}_\e$ is defined in the space $L_2(x_0,x_1)$ by the
same scheme as one used when defining this operator in the space
$L_2(\mathbb{R})$. By analogy with the proof of
Lemma~\ref{lm2.1} one can make sure that the operator
$\mathcal{H}_\e^{(0)}$ is closed for all $\e$ small enough. The
operator $\mathcal{H}_\e$ is a finite-dimensional extension of
$\mathcal{H}_0^{(-)}\oplus \mathcal{H}_\e^{(0)}\oplus
\mathcal{H}_0^{(+)}$. Hence, by \cite[Ch. I, Sec. 1.2, Thm.
4]{Gl} the identity $\conspec(\mathcal{H}_\e)=\conspec\big(
\mathcal{H}_0^{(-)}\oplus \mathcal{H}_0^{(+)}
\big)\cup\conspec(\mathcal{H}_\e^{(0)})$ is valid. It is also
follows from \cite[Ch. I, Sec. 1.1, 1.2]{Gl} that
$\resspec(\mathcal{H}_\e)\subseteq\resspec(\mathcal{H}_0^{(+)}
\oplus \mathcal{H}_\e^{(0)}\oplus \mathcal{H}_0^{(+)})=
\resspec(\mathcal{H}_0^{(+)})\cup\resspec(\mathcal{H}_\e^{(0)})
\cup\resspec(\mathcal{H}_0^{(+)})$. The operators
$\mathcal{H}_0^{(\pm)}$ are self-adjoint, and thus
$\resspec(\mathcal{H}_0^{(-)})=\resspec(\mathcal{H}_0^{(+)})=
\emptyset$. In view of Lemma~\ref{lm3.3} it sufficient to show
that $\conspec(\mathcal{H}_\e^{(0)})=
\resspec(\mathcal{H}_\e^{(0)})=\emptyset$. Let us prove these
identities.

By Lemma~\ref{lm3.2} we obtain in turn that for all $\e$ small
enough the operator $\big(\I-\e
\mathcal{L}_\e(\mathcal{H}_0^{(0)}-\iu)^{-1}\big)^{-1}:
L_2(x_0,x_1)\to L_2(x_0,x_1)$ is bounded and due to
\begin{equation*}
(\mathcal{H}_\e^{(0)}-\iu)^{-1}=(\mathcal{H}_0^{(0)}-\iu)^{-1}
\big(\I-\e\mathcal{L}_\e
(\mathcal{H}_0^{(0)}-\iu)^{-1}\big)^{-1}
\end{equation*}
the resolvent $(\mathcal{H}_\e^{(0)}-\iu)^{-1}$ is a compact
operator in $L_2(x_0,x_1)$. By \cite[Ch. I\!I\!I, Sec. 6.8, Thm.
6.29]{K} it follows that the spectrum of $\mathcal{H}_\e^{(0)}$
consists of at most countably many eigenvalues of finite
multiplicity and thus $\conspec(\mathcal{H}_\e^{(0)})=
\resspec(\mathcal{H}_\e^{(0)})=\emptyset$.
\end{proof}

\section{Proof of Theorem~\ref{th2.3}
}

By $\rho(\l)$ and $\k(\l)$ we denote the multiplier and the
quasi-momentum multiplied by $-\iu$ corresponding to the
equation (\ref{1.4}):
\begin{equation}\label{3.0}
\rho(\l):=\frac{D(\l)+\sqrt{D^2(\l)-4}}{2},\quad
\k(\l):=\ln\rho(\l).
\end{equation}
The branch of the root is specified by the requirement
$|\rho(\l)|\geqslant 1$. In the case $|\rho(\l)|=1$, the
concrete choice of the branch is not important. The branch of
logarithm is specified by $\ln 1=0$.

In accordance with Floquet-Lyapunov theorem the equation
(\ref{1.4}) has a fundamental system,
\begin{equation*}
\vp_1(x,\l)=\E^{\k(\l)x}\Phi_1(x,\l),\quad
\vp_2(x,\l)=\E^{-\k(\l)x}\Phi_2(x,\l),
\end{equation*}
where $\Phi_i(\cdot,\l)$ are 1-periodic functions.  These
formulas are valid if $\rho(\l)\not=\pm 1$, as well as in the
case $\rho(\l)=\pm 1$, $|\tht_1(\l)|^2+|\tht'_2(\l)|^2\not=0$.
If $\rho(\l)=\pm 1$, $\tht_1(\l)=\tht'_2(\l)=0$, a fundamental
system of the equation (\ref{1.4}) reads as follows
\begin{equation*}
\vp_1(x,\l)=\E^{\k(\l)x}\Phi_1(x,\l),\quad
\vp_2(x,\l)=\E^{\k(\l)x}\big(x\Phi_1(x,\l)+\Phi_2(x,\l)\big),
\end{equation*}
where $\Phi_i(\cdot,\l)$ are 1-periodic functions.

We indicate by $W[f_1,f_2]$ the Wronskian of the functions
$f_1=f_1(x)$, $f_2=f_2(x)$.

\begin{lemma}\label{it1lm1.1a}
For any $\a_1,\a_2\in\mathbb{R}$ the functions $\tht_i(x,\l)$,
$i=1,2$, are holomorphic w.r.t. $\l$ in the norm of
$C^2[\a_1,\a_2]$. The function $D(\l)$ is holomorphic. The
branches of the function $\rho(\l)$ are holomorphic everywhere
in the complex plane except the edges of non-degenerate lacunas
in the spectrum of $\mathcal{H}_0$, those are branching points
for this function. The identity
\begin{equation}\label{3.3}
W[\tht_1,\tht_2](x)=\frac{1}{p(x)}
\end{equation}
is valid.
\end{lemma}

\begin{proof}
The identity (\ref{3.3}) follows from the initial conditions
(\ref{1.3a}) and the Liouville formula for the Wronskian.

In view of $\mathcal{P}(\l,\a)$ we indicate an integral operator
\begin{equation}\label{3.2}
\big(\mathcal{P}(\l,\a)f\big)(x):=\int\limits_{\a}^x \big(
\tht_1(x,\l)\tht_2(t,\l)-\tht_1(t,\l)\tht_2(x,\l)\big)f(t)\di t.
\end{equation}
By (\ref{3.3}) this operator determines the solution to a Cauchy
problem
\begin{equation*}
\left(-\frac{d}{dx}p\frac{d}{dx}+q\right)v=f,\quad x\in
\mathbb{R}, \qquad v(\a)=\frac{dv}{dx}(\a)=0.
\end{equation*}
It is easy to make sure that for all $\l\in\mathbb{C}$ a linear
operator $\mathcal{P}: C[\a_1,\a_2]\to C^2[\a_1,\a_2]$ is
bounded for any $\a_1\leqslant\a\leqslant\a_2$. One can check
that the functions $\tht_i(x,\l)$ are solutions to the equation
\begin{equation}\label{2.25}
\big(\I-\l
\mathcal{P}(0,0)\big)\tht_i(\cdot,\l)=\tht_i(\cdot,0),
\end{equation}
which we regard as one in $C[\a_1,\a_2]$. This is the Volterra
equation. Hence the operator $\I-\l \mathcal{P}(0,0)$ is
boundedly invertible for all $\l\in\mathbb{C}$ (see, for
instance, \cite[Ch. X\!I, Sec. 3.3]{KF}). By \cite[Ch. X\!I,
Sec. 4, Proposition 4.5]{SP} it follows that the operator
$\big(\I-\l \mathcal{P}(0,0)\big)^{-1}$ is boundedly holomorphic
w.r.t. $\l\in\mathbb{C}$ as an operator in $C[\a_1,\a_2]$.
Therefore, the functions $\tht_i(x,\l)$ are holomorphic w.r.t.
$\l\in\mathbb{C}$ in the norm of $C[\a_1,\a_2]$. Since by
(\ref{2.25}) the relations
\begin{equation*}
\tht_i(x,\l)=\tht_i(x,0)+\l
\big(\mathcal{P}(0,0)\tht_i(\cdot,\l)\big)(x,\l),
\end{equation*}
hold true, we infer that the functions $\tht_i(x,\l)$ are
holomorphic w.r.t. $\l\in\mathbb{C}$ in the norm of
$C^2[\a_1,\a_2]$ as well. The holomorphy of the function $D(\l)$
is implied by one of $\tht_i(x,\l)$. The identities $D(\l)^2=4$
hold true only for $\l=\mu_n^\pm$ (see \cite[Ch. 2, Sec. 2.3,
Thm. 2.3.1]{E}), this is why only these points can be branching
points of the function $\rho$. Item c) of the proof of
Theorem~2.3.1 in \cite[Ch. 2, Sec. 2.3]{E} implies that if
$\mu_n^\pm$ is an edge of a non-degenerate lacuna, it follows
that $(D^2-4)'|_{\l=\mu_n^\pm}\not=0$. If a lacuna degenerates
($\mu_n^-=\mu_n^+$), the relations $(D^2-4)'|_{\l=\mu_n^\pm}=0$,
$(D^2-4)''|_{\l=\mu_n^\pm}\not=0$ hold true. This fact implies
the statement of the lemma on $\rho$.
\end{proof}

\begin{proof}[Proof of Theorem~\ref{th2.3}] Suppose that
$\l\in\conspec(\mathcal{H}_0)\cap K$ is an eigenvalue of the
operator $\mathcal{H}_\e$ for some $\e$ small enough. An
associated eigenfunction satisfies the equation
\begin{equation}\label{5.1}
\left(-\frac{d}{dx}p\frac{d}{dx}+q-\l-
\e\mathcal{L}_\e\right)\psi=0.
\end{equation}
For $x\not\in Q$ this equation coincides with the equation in
(\ref{1.4}). Hence, for $x$ lying to the left w.r.t. the set
$Q$, the function $\psi$ reads as follows,
\begin{equation*}
\psi(x)=c_1\vp_1(x,\l)+c_2\vp_2(x,\l),
\end{equation*}
where $c_i$ are constants. The similar behaviour is valid for
$x$ lying to the right w.r.t. $Q$. Since
$\l\in\conspec(\mathcal{H}_0)$, due to Item~(v) of Theorem~2.3.1
in \cite[Ch. 2, Sec. 2.3]{E} the identity $|\rho(\l)|=1$ is
valid, and hence $\RE\k(\l)=0$. It follows that the functions
$\vp_i(x,\l)$ are not square integrable at infinity. The
function $\psi$ is thus an element of $L_2(\mathbb{R})$ only if
\begin{equation}\label{8.1}
\psi\equiv0,\quad x\not\in Q.
\end{equation}
Since $\psi\in\H^2(\mathbb{R})\subset C^1(\mathbb{R})$, we have
\begin{equation}\label{5.3}
\psi(x_0)=\psi'(x_0)=0.
\end{equation}
We remind that the points $x_0$, $x_1$ are so that
$\overline{Q}\subset(x_0,x_1)$. By the definition (\ref{3.2}) of
$\mathcal{P}$ the initial problem (\ref{5.1}), (\ref{5.3}) is
equivalent to an integral equation
\begin{equation}\label{5.4}
\big(\I-\e \mathcal{L}_\e \mathcal{P}(\l,x_0)\big)\psi=0.
\end{equation}
The integral operator $\mathcal{P}$ is a linear bounded operator
from $L_2(Q)$ into $\H^2(Q)$. Moreover, by Lemma~\ref{it1lm1.1a}
it is bounded uniformly in  $\l\in K$. This fact and the uniform
in $\e$ boundedness of $\mathcal{L}_\e$ yield that the operator
$\mathcal{L}_\e \mathcal{P}(\l,x_0): L_2(Q)\to L_2(Q)$ is
bounded uniformly in $\e$ and $\l\in K$. Thus, for $\e$ small
enough and $\l\in K$ the operator $\big(\I-\e \mathcal{L}_\e
\mathcal{P}(\l,x_0)\big)$ is boundedly invertible, and the
equation (\ref{5.4}) therefore has the trivial solution only.
Hence, $\psi\equiv0$ for $x\in Q$. In view of (\ref{8.1}) it
follows that the function $\psi$ is identically zero. It
contradicts to the assumption that $\psi$ is an eigenfunction.
\end{proof}

\section{Proof of Theorems~\ref{th2.33},~\ref{th2.32}}

Consider the equation
\begin{equation}\label{7.4}
\left(-\frac{d}{dx}p\frac{d}{dx}+q-\e \mathcal{L}_\e-
\l\right)u=f,\quad x\in\mathbb{R},
\end{equation}
where $f\in L_2(\mathbb{R};(x_0,x_1))$,
$L_2(\mathbb{R};(x_0,x_1))$ is a subset of the functions in
$L_2(\mathbb{R})$ having supports inside $[x_0,x_1]$. We are
looking for the solutions to this equation satisfying the
conditions
\begin{equation}\label{7.3}
\begin{aligned}
&u(x,\l)=\E^{-\k(\l)x}\Phi_+(x,\l),\quad x\geqslant x_1,
\\
&u(x,\l)=\E^{\k(\l)x}\Phi_-(x,\l),\hphantom{^-}\quad x\leqslant
x_0,
\end{aligned}
\end{equation}
where $\Phi_\pm$ are 1-periodic in $x$ functions, and the branch
of the logarithm in the definition of $\k(\l)$ is specified by
the relation $\ln 1=0$. Here the branch of the function $\rho$
is not specified yet. We will study the dependence of the
solution to (\ref{7.4}), (\ref{7.3}) of $\l$, which allows us to
prove Theorems~\ref{th2.33},~\ref{th2.32}. In order to solve the
problem (\ref{7.4}), (\ref{7.3}) we employ the scheme suggested
in \cite[Ch. X\!I\!V, Sec. 4]{SP}.

We set
\begin{align*}
G(x,t,\l):=&\frac{1}{\rho(\l)-\rho^{-1}(\l)} \big(
\tht_2(\l)\tht_1(t,\l)\tht_1(x,\l)-(\rho(\l)-\tht'_2(\l))
\tht_2(t,\l)\tht_1(x,\l)
\\
&+(\rho(\l)-\tht_1(\l))\tht_1(t,\l)\tht_2(x,\l)-
\tht_1'(\l)\tht_2(t,\l)\tht_2(x,\l)\big), \quad t\geqslant x
\\
G(x,t,\l):=&\frac{1}{\rho(\l)-\rho^{-1}(\l)}
\big(\tht_2(\l)\tht_1(t,\l)\tht_2(x,\l)-(\rho(\l)-\tht'_2(\l))
\tht_1(t,\l)\tht_2(x,\l)
\\
&+(\rho(\l)-\tht_1(\l))\tht_2(t,\l)\tht_1(x,\l)-
\tht_1'(\l)\tht_2(t,\l)\tht_2(x,\l)\big), \quad x\geqslant t.
\end{align*}
The function $G(x,t,\l)$ is well-defined for all
$\l\in\mathbb{C}$ except the edges of the non-degenerate lacunas
in the spectrum of $\mathcal{H}_0$. Indeed, if
$\l\not=\mu_n^\pm$, it follows that $\rho^2(\l)\not=1$
(\cite[Ch. 2, Sec. 2.3, Thm. 2.3.1]{E}), and thus
$\rho(\l)\not=\rho^{-1}(\l)$. If lacuna degenerates
($\mu_n^-=\mu_n^+$), by Item c) of the proof of Theorem~2.3.1 in
\cite[Ch. 2, Sec. 2.3]{E} the identities
\begin{equation}\label{5.9a}
\begin{aligned}
&D(\l)=(-1)^n\big(2-\g(\l-\mu)^2\big)+\Odr\big(|\l-\mu|^3\big),
\\
&\rho(\l)=(-1)^n\big(1\pm\iu\sqrt{\g}(\l-\mu)\big)+
\Odr\big(|\l-\mu|^2\big),
\\
&\tht_1(\l)=(-1)^n+\overset{\,\textbf{.}}{\tht}_1(\mu)(\l-\mu)+
\Odr\big(|\l-\mu|^2\big),
\\
&\tht_2(\l)=\overset{\,\textbf{.}}{\tht}_2(\mu)(\l-\mu)+
\Odr\big(|\l-\mu|^2\big),
\\
&\tht'_1(\l)=\overset{\,\textbf{.}}{\tht'_1}(\mu)(\l-\mu)+
\Odr\big(|\l-\mu|^2\big),
\end{aligned}
\end{equation}
are valid, where $\l\to\mu:=\mu_n^-=\mu_n^+$, and $\g>0$ is a
constant. The sign ''$\pm$'' in the identity for $\rho(\l)$
corresponds to the different branches of this function. These
identities imply that there exists a finite limit of the
function $G$ as $\l\to\mu$, which we regard as a definition of
this function at $\l=\mu$.

On the functions $f\in L_2(\mathbb{R};(x_0,x_1))$ we introduce
the operator $\mathcal{G}(\l)$ with the kernel $G(x,t,\l)$:
\begin{equation*}
\big(\mathcal{G}(\l)f\big)(x,\l):=\int\limits_{\mathbb{R}}
G(x,t,\l)f(t)\di t.
\end{equation*}
It is clear that it is bounded as an operator from
$L_2(x_0,x_1)$ into $\H^2(x_0,x_1)$ for all values of $\l$ not
coinciding with the edges of non-degenerate lacunas in the
spectrum of $\mathcal{H}_0$.

Bearing in mind the definition of $\rho(\l)$ and $D(\l)$ by
direct calculations we check that the function $G$ is the Green
function for the equation
\begin{equation}\label{5.10}
\left(-\frac{d}{dx}p\frac{d}{dx}+q-\l\right)v=f,\quad x\in
\mathbb{R},
\end{equation}
and $v:=\mathcal{G}(\l)f$, where $f\in
L_2(\mathbb{R};(x_0,x_1))$, is a solution to this equation.

Taking into account the initial conditions  (\ref{1.3a}) and the
periodicity of $p$ and $q$, one can check easily that
\begin{equation}\label{1.11}
\tht_i(x+1,\l)=\tht_i(1,\l)\tht_1(x,\l)+\tht'_i(1,\l)\tht_2(x,\l),\quad
i=1,2.
\end{equation}
Employing this relation, the identities (\ref{3.3}) and the
formula
\begin{equation}\label{7.11}
\rho(\l)+\frac{1}{\rho(\l)}=\tht_1(\l)+\tht'_2(\l),
\end{equation}
following from (\ref{3.0}), it is easy to make sure that the
function $v$ obeys the identities
\begin{equation}\label{6.5}
v(x+1,\l)=\frac{1}{\rho(\l)}v(x,\l),\quad x\geqslant x_1,\quad
v(x-1,\l)=\frac{1}{\rho(\l)}v(x,\l),\quad x\leqslant x_0.
\end{equation}
It implies that the function $v$ satisfies the conditions
(\ref{7.3}).

Let $\z=\z(x)$ be an infinitely differentiable cut-off function
vanishing for $x\not\in[x_0,x_1]$ and equalling one in a
neighbourhood of the segment $\overline{Q}$, and $g\in
L_2(\mathbb{R};(x_0,x_1))$ be a function. We denote
\begin{equation}\label{7.5}
v:=\mathcal{G}(\l)g,\quad
w_\e:=(\mathcal{H}_\e^{(0)}-\iu)^{-1}\mathcal{L}_\e v.
\end{equation}
where, we remind, $\mathcal{H}_\e^{(0)}$ is an operator
introduced in the proof of Theorem~\ref{th2.1}. This operator is
bounded as one from the subspace of the functions in
$\H^2(x_0,x_1)$ vanishing at $x_0$, $x_1$, into $L_2(x_0,x_1)$.
As it was established in the proof of Theorem~\ref{th2.1}, the
operator $(\mathcal{H}_\e^{(0)}-\iu)^{-1}$ exists for all $\e$
small enough, and Banach theorem on inverse operator \cite[Ch.
I\!V, Sec. 5.4, Thm. 3]{KF} thus implies that for  $\e$ small
enough the operator $(\mathcal{H}_\e^{(0)}-\iu)^{-1}:
L_2(x_0,x_1)\to\H^2(x_0,x_1)$ is bounded. It is clear that it is
bounded uniformly in $\e$. Thus, an operator
$(\mathcal{H}_\e^{(0)}-\iu)^{-1}\mathcal{L}_\e$ is bounded
uniformly in $\e$ as one in $\H^2(x_0,x_1)$.

We construct the solution to the problem (\ref{7.4}),
(\ref{7.3}) as follows
\begin{equation}\label{7.6}
u(x,\l):=v(x,\l)+\e\z(x)w_\e(x,\l).
\end{equation}
This functions satisfies the conditions (\ref{7.3}). We
substitute it into the left hand side of  (\ref{7.4}) to obtain
\begin{align*}
&\left(-\frac{d}{dx}p\frac{d}{dx}+q-\e \mathcal{L}_\e-
\l\right)(v+\e\z w_\e)=g-\e \mathcal{L}_\e v+\e
\mathcal{T}_\e(\l)g
\\
&\hphantom{\big(-\frac{d}{dx}p\frac{d}{dx}}+ \e\z
\left(-\frac{d}{dx}p\frac{d}{dx}+q-\e
\mathcal{L}_\e-\iu\right)w_\e=g+\mathcal{T}_\e(\l)g,
\\
&\mathcal{T}_\e(\l)g:=-\frac{d}{dx}pw_\e\frac{d\z}{dx}-
p\frac{d\z}{dx}\frac{dw_\e}{dx}+(\iu-\l)\z w_\e.
\end{align*}
Here we have also employed the relation $\mathcal{L}_\e \z
w_\e=\mathcal{L}_\e w_\e=\z\mathcal{L}_\e w_\e$ which follows
from the identity $\z\equiv1$, $x\in\overline{Q}$. Thus, the
function $u$ defined by (\ref{7.6}) is a solution to (\ref{7.4})
if
\begin{equation}\label{7.7}
g+\e \mathcal{T}_\e(\l)g=f.
\end{equation}

\begin{lemma}\label{lm7.1}
The problem (\ref{7.4}), (\ref{7.3}) is equivalent to the
equation (\ref{7.7}) for all $\l\in \mathbb{C}$ not coinciding
with the edges of the non-degenerate lacunas in the spectrum of
$\mathcal{H}_0$. Namely, for each solution of (\ref{7.7}) there
exists the unique solution to the problem (\ref{7.4}),
(\ref{7.3}) defined by (\ref{7.5}), (\ref{7.6}). For each
solution $u$ of the problem (\ref{7.4}), (\ref{7.3}) there
exists  the unique function $g$ satisfying the equation
(\ref{7.7}) and related to $u$ by (\ref{7.5}), (\ref{7.6}).
\end{lemma}
\begin{proof}
If $g$ is a solution to (\ref{7.7}), as it was shown above, it
follows that the function $u$ introduced by  (\ref{7.5}),
(\ref{7.6}) is a solution to (\ref{7.4}), (\ref{7.3}).

Let $u$ be a solution to (\ref{7.4}), (\ref{7.3}). We define the
functions $v$, $w_\e$ and $g$ as
\begin{equation*}
w_\e:=(\mathcal{H}_0^{(0)}-\iu)^{-1} \mathcal{L}_\e u,\quad
v:=u-\e\z w_\e,\quad
g:=\left(-\frac{d}{dx}p\frac{d}{dx}+q-\l\right)v.
\end{equation*}
The function $v$ satisfies the relations (\ref{7.3}) and hence
the former of the formulas  (\ref{7.5}) holds true. The identity
(\ref{7.6}) is obviously to be valid. Since
\begin{equation*}
(\mathcal{H}_0^{(0)}-\e \mathcal{L}_\e-\iu)w_\e=\mathcal{L}_\e
u-\e \mathcal{L}_\e w_\e=\mathcal{L}_\e (u-\e\z
w_\e)=\mathcal{L}_\e v,
\end{equation*}
the latter of the formulas (\ref{7.5}) holds true as well. The
definition of the functions $g$ and $v$ and the equation
(\ref{7.4}) imply that
\begin{align*}
g&=\left(-\frac{d}{dx}p\frac{d}{dx}+q-\l\right)(u-\e\z w)=f+\e
\mathcal{L}_\e u-\e
\mathcal{T}_\e(\l)g-\e\z(\mathcal{H}_0^{(0)}-\iu)w_\e
\\
&=f+\e \mathcal{L}_\e u-\e \mathcal{T}_\e(\l)g-\e\z
\mathcal{L}_\e u=f-\e \mathcal{T}_\e(\l)g,
\end{align*}
which yields the equation (\ref{7.7}).
\end{proof}
The properties of $\mathcal{G}$ and
$(\mathcal{H}_\e^{(0)}-\iu)^{-1}\mathcal{L}_\e$ established
above allow us to claim that the mapping $g\mapsto w_\e$ is a
linear operator from $L_2(x_0,x_1)$ into $\H^1(x_0,x_1)$ bounded
uniformly in $\e$. Therefore, the operator
$\mathcal{T}_\e(\l):L_2(x_0,x_1)\to\H^1(x_0,x_1)$ is bounded
uniformly in $\e$. Moreover, it follows that the operator
$\mathcal{T}_\e$ is compact as one in $L_2(x_0,x_1)$.

\begin{proof}[Proof of Theorem~\ref{th2.32}]
For any $\d>0$ we denote
$K_\d:=K\setminus\bigcup\limits_n\big(B_\d(\mu_n^-)\cup
B_\d(\mu_n^+)\big)$, where the union is taken over the edges of
the non-degenerate lacunas in the spectrum of $\mathcal{H}_0$.
For $\l\in K_\d$ we specify the branch of the root in the
definition (\ref{3.0}) of the function $\rho$ by the condition
$|\rho(\l)|\geqslant 1$; in the case $|\rho(\l)|=1$ the exact
choice of the branch is inessential. The operator
$\mathcal{G}(\l)$ is piecewise continuous w.r.t. $\l\in K_\d$.
Therefore, the operator $\mathcal{T}_\e(\l)$ is bounded
uniformly in $\e$ and $\l\in K_\d$ as an operator in
$L_2(x_0,x_1)$. By this fact we infer that for  $\e$ small
enough the operator $\I+\e \mathcal{T}_\e(\l)$ is boundedly
invertible for all $\l\in K_\d$, and thus the equation
\begin{equation}\label{5.10b}
g+\e \mathcal{T}_\e(\l)g=0
\end{equation}
has no nontrivial solutions for $\l\in K_\d$ and $\e$ small
enough. By Lemma~\ref{lm7.1} it implies that the problem
(\ref{5.1}), (\ref{7.3}) has no nontrivial solutions. Since $\RE
\k(\l)\geqslant 0$, $\l\in K_\d$, by the choice of the branch of
$\rho$, the equation (\ref{5.1}) has no nontrivial solutions in
the space $\H^2(\mathbb{R})$. Therefore, the operator
$\mathcal{H}_\e$ has no eigenvalues in $K_\d$ if $\e$ is small
enough. This fact and the arbitrary choice of  $\d$ complete the
proof.
\end{proof}

\begin{lemma}\label{lm7.2}
The point spectrum of $\mathcal{H}_\e$ consist of at most
countably many eigenvalues. The edges of non-degenerate lacunas
in the spectrum of $\mathcal{H}_0$ are the only possible finite
accumulation points for these eigenvalues. Each eigenvalue not
coinciding with such edge is of finite multiplicity.
\end{lemma}

\begin{proof}
For any $\d>0$ we indicate $M:=\{\l\in \mathbb{C}:
\IM\l\geqslant
0\}\setminus\bigcup\limits_n\big(B_\d(\mu_n^-)\cup
B_\d(\mu_n^+)\big)$. Here the union is taken over $n$
corresponding to the edges of the non-degenerate lacunas in the
spectrum of $\mathcal{H}_0$. The statement of
Lemma~\ref{it1lm1.1a} on $\rho$ implies that the set $M$ can be
covered as follows, $M\subset
\widetilde{M}:=\bigcup\limits_{j=1}^\infty M_j$, where $M_j$ are
simply connected domains such that  $M_j\setminus
S_{\pi/4}(\mu_0^+-1)\not=\emptyset$, and for each of them one
can choose the branch of the function $\rho$ to be holomorphic
w.r.t. $\l\in M_j$, and to obey the estimate
$|\rho(\l)|\geqslant 1$ for $\IM\l\geqslant 0$.
Lemma~\ref{it1lm1.1a} and (\ref{5.9a}) yield that the operator
$\mathcal{G}(\l): L_2(x_0,x_1)\to\H^2(x_0,x_1)$ is holomorphic
w.r.t. $\l\in M_j$. This fact and the properties of the operator
$(\mathcal{H}_\e^{(0)}-\iu)^{-1}\mathcal{L}_\e$ established
above follow that the operator $\mathcal{T}_\e(\l)$ is
holomorphic in $\l\in M_j$ as an operator in $L_2(x_0,x_1)$.

Let $\l\in M$ be an eigenvalue of the operator $\mathcal{H}_\e$,
then $\l\in M_j$ for some $j$. The associated eigenfunction is a
solution to the problem (\ref{5.1}), (\ref{7.3}), where $\k$ is
defined via the branch of $\rho$ corresponding to $M_j$. Due to
Lemma~\ref{lm7.1} it means that the corresponding equation
(\ref{5.10b}) has a nontrivial solution. The compactness of the
operator $\mathcal{T}_\e(\l)$ implies that the number of such
linear independent solutions is finite and $\l$ is thus an
eigenvalue of finite multiplicity. Since $M_j\setminus
S_{\pi/4}(\mu_0^+)\not=\emptyset$, in accordance with
Theorem~\ref{th2.2} there exists a point $\l_*\in
M_j\setminus\spec(\mathcal{H}_\e)$. Therefore, the equation
(\ref{5.10b}) has no nontrivial solution for $\l=\l_*$ that
together with the compactness of $\mathcal{T}_\e(\l_*)$ implies
the bounded invertibility of $\I+\e \mathcal{T}_\e(\l_*)$. This
fact and the above established holomorphy of $\mathcal{T}_\e$ in
$\l\in M_j$ allow us to employ Theorem~7.1 in \cite[Ch. X\!V,
Sec. 7]{SP} and to conclude that the operator $(\I+\e
\mathcal{T}_\e)^{-1}$ is meromorphic in $\l\in M_j$ and has at
most countably many poles in $M_j$ those can accumulate at the
boundary points of  $M_j$ only. Moreover, the poles of $(\I+\e
\mathcal{T}_\e)^{-1}$ are values of $\l$ for which the equation
(\ref{5.10b}) has a nontrivial solution. Therefore, the
eigenvalues of $\mathcal{H}_\e$ lying in $M_j\cap\{\l\in
\mathbb{C}: \IM\l\geqslant 0 \}$ are poles of $(\I+\e
\mathcal{T}_\e)^{-1}$ corresponding to $M_j$. Thus, the operator
$\mathcal{H}_\e$ has at most countably many eigenvalues in $M$.
The points of $M$ can not be accumulation points for these
eigenvalues since each such point is inner for one of the sets
$M_j$. In the same way one can prove that the set obtained from
$M$ by mirror symmetry w.r.t. real axis contains at most
finitely many eigenvalues of $\mathcal{H}_\e$ of finite
multiplicity those have no accumulation points inside this set.
The number $\d$ being arbitrary completes the proof.
\end{proof}

In view of this lemma it remains to prove that edges of the
non-degenerate lacunas in the spectrum of $\mathcal{H}_0$ are
not the accumulation points for the eigenvalues of
$\mathcal{H}_\e$. We should also show that in the case such an
edge is an eigenvalue of $\mathcal{H}_\e$ it is of finite
multiplicity. We will prove these facts on the basis of the
equation similar to (\ref{7.7}). We can not employ exactly this
equation since the function $\rho$ has branching points at
$\mu_n^\pm$, and hence the operator $\mathcal{T}_\e$ is not
holomorphic at these points.

First we prove an auxiliary statement.

\begin{lemma}\label{lm7.3}
Let  $\mu_n^\pm$ be an edge of a non-degenerate lacuna in the
spectrum of $\mathcal{H}_0$. Then
\begin{enumerate}\def\theenumi{\arabic{enumi}}

\item\label{it4lm1.1}  At least one of the numbers
$\tht'_1(\mu_n^+)$ and $\tht_2(\mu_n^+)$ (respectively,
$\tht'_1(\mu_n^-)$ and $\tht_2(\mu_n^-)$) is non zero and the
inequality $\tht'_1(\mu_n^\pm)\tht_2(\mu_n^\pm)\leqslant 0$
holds true.

\item\label{it2alm1.1}
$D(\mu_n^\pm)=2(-1)^n$,
$\mp(-1)^n\overset{\,\textbf{.}}{D}(\mu_n^\pm)>0$.

\item\label{it2clm1.1}
The numbers $\overset{\,\textbf{.}}{D}(\mu_n^\pm)$ are given by
\begin{align*}
& \overset{\,\textbf{.}}{D}(\mu_n^\pm)=\int\limits_0^1
\Big(\tht'_1(\mu_n^\pm)\tht_2^2(x,\mu_n^\pm)+
(\tht_1(\mu_n^\pm)-\tht'_2(\mu_n^\pm))\tht_1(x,\mu_n^\pm)
\tht_2(x,\mu_n^\pm)
-
\\
&\hphantom{\overset{\,\textbf{.}}{D}(\mu_n^\pm)=\int\limits_0^1
\Big( }-\tht_2(\mu_n^\pm)\tht_1^2(x,\mu_n^\pm)\Big)^2\di x.
\end{align*}
Moreover,
\begin{align*}
&\overset{\,\textbf{.}}{D}(\mu_n^\pm)=-
\frac{1}{4\tht_2(\mu_n^\pm)} \int\limits_0^1
\Big(2\tht_2(\mu_n^\pm)\tht_1(x,\mu_n^\pm)+
(\tht_1(\mu_n^\pm)-\tht'_2(\mu_n^\pm))\tht_2(x,\mu_n^\pm)\Big)^2\di
x,
\end{align*}
if $\tht_2(\mu_n^\pm)\not=0$, and
\begin{equation*}
\overset{\,\textbf{.}}{D}(\mu_n^\pm)=
\frac{1}{4\tht'_1(\mu_n^\pm)} \int\limits_0^1
\Big(2\tht'_1(\mu_n^\pm)\tht_2(x,\mu_n^\pm)+
(\tht_1(\mu_n^\pm)-\tht'_2(\mu_n^\pm))\tht_1(x,\mu_n^\pm)\Big)^2\di
x,
\end{equation*}
if $\tht'_1(\mu_n^\pm)\not=0$.
\end{enumerate}
\end{lemma}

The lemma follows from Theorem~2.3.1 in \cite[Ch. 2, Sec.
2.3]{E}, and the formula (2.3.7) and the relation
$D^2(\l)=4+(\tht_1(\l)-\tht'_2(\l))^2+4\tht'_1(\l)\tht_2(\l)$
established in the proof of this theorem.

Let  $\mu_n^\pm$ be an edge of a non-degenerate lacuna in the
spectrum of $\mathcal{H}_0$.  In a small neighbourhood of
$\mu_n^\pm$ we introduce a new complex parameter by the rule
$\l:=\mu_n^\pm\mp k^2$. We denote
\begin{equation}\label{5.15c}
\rho_n^\pm=\rho_n^\pm(k):=\frac{D(\mu_n^\pm\mp k^2)+(-1)^n
\sqrt{D^2(\mu_n^\pm\mp k^2)-4}}{2},
\end{equation}
where the branch of the root is specified by the condition
$\sqrt{1}=1$, if $\arg k\in[0,\pi)$ and $\sqrt{1}=-1$, if $\arg
k\in[\pi,2\pi)$. By Item~\ref{it2alm1.1} of Lemma~\ref{lm7.3}
the identities
\begin{equation*}
D(\mu_n^\pm\mp
k^2)=(-1)^n(2+|\overset{\,\textbf{.}}{D}(\mu_n^\pm)|
k^2)+\Odr(|k|^4),\quad D^2(\mu_n^\pm\mp
k^2)-4=4|\overset{\,\textbf{.}}{D}(\mu_n^\pm)| k^2+\Odr(|k|^4)
\end{equation*}
are valid for $k$ small enough. Thus, the function $\rho_n^\pm$
is holomorphic in $k$ and its Taylor expansion reads as follows,
\begin{equation}\label{1.9}
\rho_n^\pm(k)=(-1)^n\left(1+\sqrt{|\overset{\,\textbf{.}}{D}(\mu_n^\pm)|}k+\frac{1}{2}|\overset{\,\textbf{.}}{D}(\mu_n^\pm)|
k^2\right)+\Odr(|k|^3).
\end{equation}
We set
\begin{align*}
&\k_n=\k_n^\pm(k):=\ln\rho_n^\pm(k),\hphantom{\pi\iu\ \ \ .}
\quad\text{if $n$ is even},
\\
&\k_n=\k_n^\pm(k):=\ln\rho_n^\pm(k)-\pi\iu,\quad\text{if $n$ is
odd},
\end{align*}
where the branch of the logarithm is specified by  $\ln 1=0$.
For $k$ small enough the function $\k_n^\pm(k)$ is holomorphic
w.r.t. $k$ and the identity
\begin{equation}\label{1.10}
\k_n^\pm(k)=\sqrt{|\overset{\,\textbf{.}}{D}(\mu_n^\pm)|}k+\Odr(|k|^2)
\end{equation}
holds true.

Consider the equation (\ref{7.4}) for $\l=\mu_n^\pm\mp k^2$. We
are looking for the solution to this equation satisfying the
conditions
\begin{equation}\label{7.10}
\begin{aligned}
&u(x,\l)=\E^{-\k_n^\pm(k)x}\Phi^\pm_{n,+}(x,k),\quad x\geqslant
x_1,
\\
&u(x,\l)=\E^{\k_n^\pm(k)x}\Phi^\pm_{n,-}(x,k),\hphantom{^-}\quad
x\leqslant x_0,
\end{aligned}
\end{equation}
where $\Phi^\pm_{n,\pm}$ are 1-periodic w.r.t. $x$ functions. In
order to solve this problem we again employ the scheme borrowed
from \cite[Ch. X\!I\!V, Sec. 4]{SP}. The main difference is that
the analogue of the function $v$ in (\ref{7.5}) is defined in a
more complicated way that allows us to avoid singularities at
the point $k=0$ for an analogue of the operator $\mathcal{T}_\e$
in (\ref{7.7}). To define the analogue of the function $v$ we
first introduce additional notations.

We denote $\tau_n^\pm:=\pm(-1)^n$, and
\begin{equation}
\begin{aligned}
&\vp_{n,1}^\pm(x,k):=\sqrt{\tau_n^\pm\tht_2(\l)}
\left(\tht_1(x,\l)+
\frac{\rho_n^\pm(k)-\tht_1(\l)}{\tht_2(\l)}\tht_2(x,\l)\right),
\\
&\vp_{n,2}^\pm(x,k):=\sqrt{\tau_n^\pm\tht_2(\l)}
\left(\tht_1(x,\l)+
\frac{\big(\rho_n^\pm(k)\big)^{-1}-\tht_1(\l)}{\tht_2(\l)}\tht_2(x,\l)\right),
\end{aligned}\label{5.15a}
\end{equation}
if $\tht_2(\mu_n^\pm)\not=0$, and
\begin{equation}
\begin{aligned}
&\vp_{n,1}^\pm(x,k):=\sqrt{-\tau_n^\pm\tht'_1(\l)}
\left(\frac{\rho_n^\pm(k)-\tht'_2(\l)}{\tht'_1(\l)}\tht_1(x,\l)+
\tht_2(x,\l)\right),
\\
&\vp_{n,2}^\pm(x,k):=\sqrt{-\tau_n^\pm\tht'_1(\l)}
\left(\frac{\big(\rho_n^\pm(k)\big)^{-1}-\tht'_2(\l)}{\tht'_1(\l)}\tht_1(x,\l)+
\tht_2(x,\l)\right),
\end{aligned}\label{5.15b}
\end{equation}
if $\tht_2(\mu_n^\pm)=0$. Everywhere in these formulas the
symbol $\l$ indicates the sum $\mu_n^\pm\mp k^2$.
Item~\ref{it4lm1.1} of Lemma~\ref{lm7.3} implies that the
functions $\vp_{n,i}^\pm$ are well-defined.

By analogy with (\ref{1.11}), (\ref{7.11}), (\ref{6.5}) one can
check that the functions $\vp_{n,i}^\pm$ can be represented as
\begin{equation}\label{1.12}
\vp_{n,1}^\pm(x,k)=\E^{\k_n^\pm(k)x}\Phi_{n,1}^\pm(x,k),\quad
\vp_{n,2}^\pm(x,k)=\E^{-\k_n^\pm(k)x}\Phi_{n,2}^\pm(x,k),
\end{equation}
where $\Phi_{n,i}^\pm(x,k)$ are 1-periodic w.r.t. $x$ for even
$n$ and 1-antiperiodic for odd $n$. Lemma~\ref{it1lm1.1a} and
1-(anti)periodicity of $\Phi_{n,i}^\pm$ imply that for $\e$
small enough these functions are holomorphic w.r.t. $k$ small
enough in the norm of $C^2[\a_1,\a_2]$ for all $\a_1,\a_2\in
\mathbb{R}$. Since the functions $\tht_i(x,\mu_n^\pm)$ are real,
Items~\ref{it2alm1.1},~\ref{it2clm1.1} of Lemma~\ref{lm7.3}
imply that
\begin{align*}
&\sgn \tht_2(\mu_n^\pm)=\pm(-1)^n,\qquad \text{if\ \
$\tht_2(\mu_n^\pm)\not=0$},
\\
-&\sgn\tht'_1(\mu_n^\pm)=\pm(-1)^n,\qquad \text{if\ \
$\tht'_1(\mu_n^\pm)\not=0$}.
\end{align*}
Bearing in mind these relations, the holomorphy of the functions
$\Phi_{n,i}^\pm$, (\ref{7.11}), (\ref{3.3}), (\ref{1.9}),
(\ref{1.10}), one can make sure that
\begin{equation}\label{1.8a}
\vp_{n,i}^\pm(x,0)=\phi_{n}^\pm(x),\quad i=1,2,
\end{equation}
where $\phi_{n}^\pm$ are the eigenfunctions of (\ref{1.3})
associated with $\mu_n^\pm$  and satisfying the normalization
condition (\ref{1.3c}). The right hand side of this relation is
nonzero by Item~\ref{it4lm1.1} of Lemma~\ref{lm7.3}. The
function $\phi_{n}^\pm$ being not identically zero, there exists
a point $x_2\not\in\overline{Q}$ such that
$\phi_{n}^\pm(x_2)\not=0$. Of course, the point $x_2$ depends on
$n$ and an edge of a lacuna. Enlarging if needed the interval
$(x_0,x_1)$ we can assume that $x_2\in(x_0,x_1)$ and $\z\equiv
1$ in a neighbourhood of the point $x_2$.

On the functions $f\in L_2(\mathbb{R};(x_0,x_1))$ we define the
operators
\begin{align*}
&\mathcal{G}_{n,+}^\pm(k)f:=\mathcal{P}(\l,+\infty)f
-\frac{\big(\mathcal{P}(\l,+\infty)f\big)(x_2,k)}
{\vp_{n,2}^\pm(x_2,k)} \vp_{n,2}^\pm(\cdot,k),
\\
&\mathcal{G}_{n,-}^\pm(k)f:=\mathcal{P}(\l,-\infty)f
-\frac{\big(\mathcal{P}(\l,-\infty)f\big)(x_2,k)}
{\vp_{n,1}^\pm(x_2,k)} \vp_{n,1}^\pm(\cdot,k),
\end{align*}
where $\l=\mu_n^\pm\mp k^2$, and $\mathcal{P}$, we remind, is
the operator in (\ref{3.2}). For the brevity till the end of the
section we will omit the index ''$\pm$'' in the notations.

By Lemma~\ref{it1lm1.1a}, the identity (\ref{1.8a}) and the
assumption $\phi_{n}^\pm(x_2)\not=0$ the operators
$\mathcal{G}_{n,\pm}^\pm$ are holomorphic w.r.t. $k$ small
enough as the operators from $L_2(\mathbb{R};(x_0,x_1))$ into
$\H^2(x_0,x_1)$. Let $g\in L_2(\mathbb{R};(x_0,x_1))$ be a
function. It is easy to check that the functions
\begin{equation*}
v_+:=\mathcal{G}_{n,+}(k)g,\quad v_-:=\mathcal{G}_{n,-}(k)g,
\end{equation*}
are solutions to the equations (\ref{5.10}) for $\l=\mu_n^\pm\mp
k^2$ vanishing at $x_2$. Moreover, the function $v_+$ satisfies
the former of the relations (\ref{7.10}), while $v_-$ does the
latter.

We introduce the function $v(x,k):=v_-(x,k)$, $x<x_2$,
$v(x,k):=v_+(x,k)$, $x>x_2$. Let $w_\e$ be a solution to the
boundary value problem
\begin{equation}\label{7.13}
\begin{gathered}
\left(-\frac{d}{dx}p\frac{d}{dx}+q-\e
\mathcal{L}_\e-\iu\right)w_\e=h,\quad x\in(x_0,x_1), \quad
w_\e=v,\quad x=x_0,x_1,
\\
h:= \left\{
\begin{aligned}
&\left(-\frac{d}{dx}p\frac{d}{dx}+q-\iu\right)v_-,& x<x_2,
\\
&\left(-\frac{d}{dx}p\frac{d}{dx}+q-\iu\right)v_+,& x>x_2.
\end{aligned}
\right.
\end{gathered}
\end{equation}
The function $h$ belongs to $L_2(x_0,x_1)$. The problem
(\ref{7.13}) is uniquely solvable in $\H^2(x_0,x_1)$ for $\e$
small enough. Indeed, a change
\begin{equation*}
w_\e=w_\e^{(0)}+w_\e^{(1)},\quad
w_\e^{(0)}(x):=\frac{v(x_1)(x-x_0)-v(x_0)(x-x_1)}{x_1-x_0},
\end{equation*}
reduces the problem (\ref{7.13}) to an equation
$(\mathcal{H}_\e^{(0)}-\iu) w_\e^{(1)}=h^{(1)}$, where
\begin{equation*}
h^{(1)}:=h-\left(-\frac{d}{dx}p\frac{d}{dx}+q-\e
\mathcal{L}_\e-\iu\right)w_\e^{(0)}.
\end{equation*}
As it was shown above this equation is uniquely solvable.
Moreover, the mapping $g\mapsto w_\e$ is a linear operator from
$L_2(\mathbb{R};(x_0,x_1))$ into $\H^2(x_0,x_1)$ being
holomorphic w.r.t. $k$ small enough.

We construct the solution to (\ref{7.4}), (\ref{7.10}) as
follows
\begin{equation}\label{7.14}
u(x,k):=\big(1-\z(x)\big)v(x,k)+\z(x)w_\e(x,k).
\end{equation}
This functions satisfies the conditions (\ref{7.10}). We
substitute it into the left hand side of (\ref{7.4}) with
$\l=\mu_n^\pm\mp k^2$ to obtain
\begin{align*}
&\left(-\frac{d}{dx}p\frac{d}{dx}+q-\e \mathcal{L}_\e-
\mu_n^\pm\pm k^2\right)\big((1-\z)v+\z
w_\e\big)=g+\mathcal{T}_{n,\e}^\pm(k)g,
\\
&\mathcal{T}_{n,\e}^\pm(k)g:=-\frac{d}{dx}p(w_\e-v)\frac{d\z}{dx}-
p\frac{d\z}{dx}\frac{d}{dx}(w_\e-v)+(\iu-\mu_n^\pm\pm k^2)\z
(w_\e-v),
\end{align*}
that leads us to the equation
\begin{equation}\label{7.15}
g+\mathcal{T}_{n,\e}^\pm(k)g=f.
\end{equation}
The operator $\mathcal{T}_{n,\e}^\pm$ is compact in
$L_2(x_0,x_1)$.

\begin{lemma}\label{lm7.4}
The problem (\ref{7.4}), (\ref{7.10}) is equivalent to the
equation (\ref{7.15}) for all $k$ small enough. Namely, for each
solution of (\ref{7.15}) there exists the unique solution to
(\ref{7.4}), (\ref{7.10}) defined by (\ref{7.14}). For each
solution $u$ of (\ref{7.4}), (\ref{7.10}) there exists  the
unique function $g$ satisfying to the equation(\ref{7.15}) and
related with $u$ by (\ref{7.14}).
\end{lemma}

The proof of the lemma is completely analogous to that of
Lemma~\ref{lm7.1}. The formulas expressing $v_\pm$, $w_\e$, $g$
via the solution $u$ of (\ref{7.4}), (\ref{7.10}) are as
follows,
\begin{align*}
&v_\pm(x,k):=u(x,k)-\z(x)U_\pm(x,k),\,\,\quad \pm x>\pm x_2,
\\
&w_\e(x,k):=u+\big(1-\z(x)\big)U_\pm(x,k),\quad \pm x>\pm x_2,
\\
&g(x,k)=\left\{
\begin{aligned}
&\left(-\frac{d}{dx}p\frac{d}{dx}+q-\mu_n^\pm\pm k^2\right)v_-,&
x<x_2,
\\
&\left(-\frac{d}{dx}p\frac{d}{dx}+q-\mu_n^\pm\pm k^2\right)v_+,&
x>x_2,
\end{aligned}
\right.
\end{align*}
where $U_\pm$ are solutions to the boundary value problems:
\begin{gather*}
\left(-\frac{d}{dx}p\frac{d}{dx}+q-\e
\mathcal{L}_\e-\iu\right)U_-=\e \mathcal{L}_\e u,\quad
x\in(x_0,x_2),
\\
U_-=0,\quad x=x_0,\quad U_-=u,\quad x=x_1,
\\
\left(-\frac{d}{dx}p\frac{d}{dx}+q-\e
\mathcal{L}_\e-\iu\right)U_+=\e \mathcal{L}_\e u,\quad
x\in(x_2,x_1),
\\
U_+=u,\quad x=x_2,\quad U_+=0,\quad x=x_2.
\end{gather*}

\begin{proof}[Proof of Theorem~\ref{th2.33}]
As it was said above, to finish the proof of
Theorem~\ref{th2.33} it remains to show that the edges of the
non-degenerate lacunas in the spectrum of  $\mathcal{H}_0$ are
not accumulation points for the eigenvalues of $\mathcal{H}_\e$,
and we  should also check that if such an edge is an eigenvalue,
it is of finite multiplicity. We will prove these facts by
analogy with the proof of Lemma~\ref{lm7.2}.

Let $\l$ be an eigenvalue of the operator $\mathcal{H}_\e$ lying
in a vicinity of an edge $\mu_n^\pm$ of a non-degenerate lacuna.
Then the corresponding eigenfunction is a nontrivial solution to
the problem
\begin{equation}\label{5.19}
\left(-\frac{d}{dx}p\frac{d}{dx}+q-\e
\mathcal{L}_\e-\mu_n^\pm\pm k^2\right)\psi=0,\quad
x\in\mathbb{R},
\end{equation}
satisfying the condition (\ref{7.10}), where
$\RE\k_n^\pm(k)\geqslant0$.  By (\ref{1.10}) the last inequality
is equivalent to $\RE k \geqslant0$. Due to Lemma~\ref{lm7.4} it
means that the corresponding equation (\ref{7.15}) has a
nontrivial solution. The compactness of the operator
$\mathcal{T}_{n,\e}^\pm(k)$ follows that the number of such
solutions is finite. Thus, if $\mu_n^\pm$ is an eigenvalue of
$\mathcal{H}_\e$, it is of finite multiplicity.

Since by Lemma~\ref{lm7.2} the number of the eigenvalues of
$\mathcal{H}_\e$ is at most countable, it follows that in any
small neighbourhood of $\mu_n^\pm$ there exists a point $\l_*$
not belonging to the spectrum of  $\mathcal{H}_\e$. Let $k_*$ be
a value of the parameter $k$ associated with $\l_*$ and $\RE
k_*>0$. Then the problem (\ref{5.19}), (\ref{7.10}) with $k=k_*$
has no nontrivial solutions that by Lemma~\ref{lm7.4} means the
absence of the nontrivial solutions to the equation (\ref{7.15})
with $k=k_*$. Thus, the operator $(\I+\mathcal{T}_\e(k_*))$ is
boundedly invertible. Together with the holomorphy of this
operator it allows us to apply Theorem~7.1 in \cite[Ch. X\!V,
Sec. 7]{SP} and to conclude that the operator
$(\I+\mathcal{T}_\e)^{-1}$ is meromorphic in $k$ and its poles
can not accumulate at the internal point of the considered
neighbourhood. Each value of the parameter $k$ corresponding to
the eigenvalue of operator $\mathcal{H}_\e$ being close to
$\mu_n^\pm$ is a pole of the operator
$(\I+\mathcal{T}_\e)^{-1}$. Thus, the eigenvalues of the
operator $\mathcal{H}_\e$ lying in the vicinity of $\mu_n^\pm$
can not accumulate at $\mu_n^\pm$.
\end{proof}

\section{Proof of Theorem~\ref{th2.31}}

In the present section we prove Theorem~\ref{th2.31}. Moreover,
we provide an example showing that under violation of the
hypothesis of this theorem the operator  $\mathcal{H}_\e$ can
have embedded eigenvalues.

By analogy with the proof of Theorem~\ref{th2.3} one can show
that if an eigenvalue $\l$ of $\mathcal{H}_\e$ is embedded, the
associated eigenfunction satisfies the identity (\ref{8.1}). In
what follows we regard this identity as proven.

\begin{proof}[Proof of Item~\ref{it1th2.4} of Theorem~\ref{th2.31}]
We argue by contradiction. Let $\l\in\conspec(\mathcal{H}_\e)$
be an eigenvalues of $\mathcal{H}_\e$, and $\psi$ is an
associated eigenfunction. Let $x_2$ be the left endpoint of the
support of $\psi$ and
\begin{equation}\label{8.2}
\psi(x)\not\equiv0,\quad x_2<x<x_2+\d,
\end{equation}
where $\d$ is a small number. Since $\psi\in
\H^2(\mathbb{R})\subset C^1(\mathbb{R})$, the initial conditions
\begin{equation}\label{6.01}
\psi(x_2)=\psi'(x_2)=0
\end{equation}
hold true. The definition (\ref{3.2}) of $\mathcal{P}$ implies
that the initial problem (\ref{5.1}), (\ref{6.01}) is equivalent
to an integral equation
\begin{equation}\label{8.3}
\big(\I-\e \mathcal{L}_\e \mathcal{P}(\l,x_2)\big)\psi=0.
\end{equation}
In view of (\ref{3.2}), Lemma~\ref{it1lm1.1a} and the estimate
(\ref{1.3e}), for any function $u\in \H^2(x_2,x_2+\d)$ we obtain
\begin{equation*}
\|\mathcal{L}_\e\mathcal{P}(\l,x_2)u\|_{L_2(x_2,x_2+\d)}\leqslant
C\d\|u\|_{L_2(x_2,x_2+\d)},
\end{equation*}
where the constant $C$ is independent of $\d$ and $\e$.
Therefore, for $\d$ small enough the operator
$\mathcal{L}_\e\mathcal{P}(\l,x_2)$ is a contraction operator in
$L_2(x_2,x_2+\d)$ and the equation (\ref{8.3}) thus has the
trivial solution only. It contradicts to (\ref{8.2}).
\end{proof}

\begin{proof}[Proof of Item~\ref{it2th2.4} of Theorem~\ref{th2.31}]
In view of Theorem~\ref{th2.3} it is sufficient to show that
there are no embedded eigenvalues of $\mathcal{H}_\e$ tending to
infinity as $\e\to0$.

Let $\l\in\conspec(\mathcal{H}_\e)$ is an eigenvalue of
$\mathcal{H}_\e$ greater than one, and $\psi$ is an associated
eigenfunction normalized in $L_2(\mathbb{R})$. This
eigenfunction satisfies the equation
\begin{equation}\label{8.4}
\left(-\frac{d}{dx}p_\e\frac{d}{dx}+q-\l-
\e\widetilde{\mathcal{L}}_\e\right)\psi=0,
\end{equation}
where $p_\e:=p+\e a_\e$. We multiply this equation by
$\overline{\psi}$ and integrate by parts bearing in mind
(\ref{8.1}). It results in
\begin{equation}\label{8.5}
\|\sqrt{p_\e}\psi'\|_{L_2(Q)}^2+(q\psi,\psi)_{L_2(Q)}+
\e(\widetilde{\mathcal{L}}_\e\psi,\psi)_{L_2(Q)}=\l.
\end{equation}
The function $p_\e$ is positive for $\e$ small enough since by
(\ref{1.3d}) and the function $a_\e$ being compactly supported
we have
\begin{equation*}
\e \max\limits_{\overline{Q}} | a_\e(x)|=\e
\max\limits_{\overline{Q}} \bigg|\int\limits_{x_0}^x a_\e'(t)\di
t\bigg| \leqslant \e |x_1-x_0| \max\limits_{\overline{Q}}
|a_\e'(t)|\xrightarrow[\e\to0]{}0.
\end{equation*}
Bearing in mind this convergence, by (\ref{8.5}), (\ref{0.1})
and the uniform in $\e$ boundedness of the operator
$\widetilde{\mathcal{L}}_\e:\H^1(Q)\to L_2(Q)$ we obtain
\begin{equation}
\|\psi'\|_{L_2(Q)}^2\leqslant C+\l+C\e
\|\psi'\|_{L_2(Q)}^2,\quad \|\psi\|_{\H^1(Q)}\leqslant
C\sqrt{\l}, \quad
\|\widetilde{\mathcal{L}}_\e\psi\|_{L_2(Q)}\leqslant
C\sqrt{\l},\label{8.7}
\end{equation}
where $C$ are constants independent of $\e$ and $\l$. We
multiply the equation (\ref{8.4}) by $\E^{\a
x}p_\e\overline{\psi}'$, $\a>0$, and integrate by parts taking
into account (\ref{8.1}):
\begin{equation}
\begin{aligned}
0=&-2\RE\int\limits_Q\E^{\a x}p_\e\overline{\psi}'(p_\e
\psi')'\di x +2\RE\int\limits_Q q\E^{\a x}p_\e
\psi\overline{\psi}'\di x+2\e \RE\int\limits_Q \E^{\a
x}p_\e\overline{\psi}' \widetilde{\mathcal{L}}_\e\psi\di x
\\
& -2\l\RE\int\limits_Q\E^{\a x}p_\e\overline{\psi}'\psi\di x
=\a\|\E^{\frac{\a\cdot}{2}}p_\e\psi'\|_{L_2(Q)}^2+\l\big((\E^{\a
x}p_\e)'\psi,\psi\big)_{L_2(Q)}
\\
&+2\RE\big( q\E^{\a
x}p_\e\overline{\psi}',\psi\big)_{L_2(Q)}+2\e \RE\big(\E^{\a
x}p_\e\widetilde{\mathcal{L}}_\e\psi,\psi'\big)_{L_2(Q)}.
\end{aligned}\label{8.8}
\end{equation}
Since $(\E^{\a x}p_\e)'=\E^{\a x}(\a p_\e+p_\e')$, by
(\ref{0.1}) and (\ref{1.3d}) we infer that there exists a
constant $\a$ independent of $\e$ and $\l$ such that
\begin{equation*}
\big(\E^{\a x}(p+\e a_\e)\big)'\geqslant 1,\quad x\in
\overline{Q}.
\end{equation*}
We substitute this estimate and (\ref{8.7}) into (\ref{8.8}) and
get
\begin{equation*}
\l\leqslant  2\left|\big( q\E^{\a x}(p+\e
a_\e)\overline{\psi}',\psi\big)_{L_2(Q)}\right|+2\e
\left|\big(\E^{\a x}(p+\e
a_\e)\widetilde{\mathcal{L}}_\e\psi,\psi'\big)_{L_2(Q)}\right|
\leqslant C\sqrt{\l}+\e C\l,
\end{equation*}
where the constant $C$ is independent of $\e$ and $\l$.
Therefore, $\l\leqslant C\sqrt{\l}$, that implies $\l\leqslant
C$, where the constant $C$ is independent of $\e$ and $\l$.
\end{proof}

In the remaining part of the section we give the example showing
that under violation of the hypothesis of this theorem the
operator  $\mathcal{H}_\e$ can have an eigenvalue embedded into
the continuous spectrum.

We set $p\equiv1$, $q\equiv0$, i.e.,
$\mathcal{H}_0=-\frac{d^2}{dx^2}$. We introduce the operator
$\mathcal{L}_\e$ as follows,
\begin{align*}
&(\mathcal{L}_\e u)(x):=2\xi_\e(x)l_\e u, && l_\e
u:=\e^{-1}\big(u'(\e^\a)-u'(0)\big), && \a\geqslant 2,
\\
&\xi_\e(x):=c_\e\chi_\e(x)\sin\left( \nu_\e|x|\right), &&
c_\e:=\left(\frac{2\pi[\nu_\e]}{\nu_\e}-\e^\a\right)^{-1}, &&
\nu_\e:=\frac{\pi}{2\e^\a},
\end{align*}
where $\chi_\e$ is the characteristic function of the interval
$\left(-\frac{2\pi[\nu_\e]}{\nu_\e},\frac{2\pi[\nu_\e]}{\nu_\e}\right)$,
$[\nu_\e]$ is the integer part of $\nu_\e$. Let us check that
the operator $\mathcal{L}_\e$ satisfies all needed requirements
with  $Q=(-2\pi,2\pi)$. It is clear that due to embedding
$\H^2(Q)\subset C^1(Q)$ the functional $l_\e$ is a linear
functional on the space $\H^2(Q)$, and $\mathcal{L}_\e$ is a
linear operator from $\H^2(Q)$ into $L_2(\mathbb{R};Q)$. By the
estimate
\begin{equation*}
|l_\e u|\leqslant \e^{-1}\int\limits_0^{\e^\a}|u''(t)|\di t
\leqslant \e^{\a/2-1}\|u''\|_{L_2(0,\e)}\leqslant
C\|u\|_{\H^2(Q)},
\end{equation*}
where the constant $C$ is independent of $\e$, we obtain the
uniform in $\e$ boundedness for the functional $l_\e$ on
$\H^2(Q)$. It yields the same property for  $\mathcal{L}_\e$.

Let us prove now that $\l_\e:=\nu_\e^2$ is an eigenvalue of the
operator $\mathcal{H}_\e$, and
\begin{equation*}
\psi_\e(x):=-\frac{\e}{\nu_\e}\int\limits_Q \sin\left(
\nu_\e|x-t|\right)\xi_\e(t)\di t
\end{equation*}
is an associated eigenfunction. Clearly,
$\psi_\e\in\Hloc^2(\mathbb{R})$. Employing the formula
\begin{equation*}
\psi_\e'(x)=-\e\int\limits_Q \cos
\left(\nu_\e(x-t)\right)\sgn(x-t)\xi_\e(t)\di t,\quad x\in Q,
\end{equation*}
by direct calculations we check that  $\psi_\e'(0)=0$,
$\psi_\e'(\e^\a)=\e$, $l_\e\psi_\e=1$. This equality and the
definition of the function $\psi_\e$ imply that this function is
a solution to (\ref{5.1}) for $\l=\l_\e$. The definition of
$\psi_\e$ also implies that $\psi_\e(x)\equiv0$ for $x\not\in
Q$, and thus $\psi_\e\in\H^2(\mathbb{R})$. Therefore, $\l_\e$ is
an eigenvalue of $\mathcal{H}_\e$. Since
$\conspec(\mathcal{H}_0)=[0,+\infty)$, by Theorem~\ref{th2.1} we
have an identity $\conspec(\mathcal{H}_\e)=[0,+\infty)$. The
eigenvalue $\l_\e$ is positive, and thus
$\l_\e\in\conspec(\mathcal{H}_\e)$.

\section{Auxiliary statements}

In the present section we prove certain auxiliary statements
needed in proof of Theorems~\ref{th2.5},~\ref{th2.6}.

Throughout the section by $\mu_n^\pm$ we mean an edge of a
non-degenerate lacuna in the spectrum of $\mathcal{H}_0$.

The identity (\ref{3.3}) and the definition (\ref{5.15a}),
(\ref{5.15b}) of the functions $\vp_{n,i}^\pm$ imply that
\begin{equation}\label{1.13a}
W[\vp_{n,1}^\pm,\vp_{n,2}^\pm]=\frac{\tau_n^\pm}{p(x)} \left(
\frac{1}{\rho_n^\pm(k)}-\rho_n^\pm(k)\right),
\end{equation}
where, we remind, the functions $\rho_n^\pm$ were introduced in
(\ref{5.15c}). By $\mathcal{G}_n^\pm(k)$ we denote an integral
operator defined in $L_2(Q)$:
\begin{gather*}
(\mathcal{G}_n^\pm(k)f)(x):=\int\limits_{Q}
G_n^\pm(x,t,k)f(t)\di t,
\\
G_n^\pm(x,t,k):=\frac{\tau_n^\pm}{
\rho_n^\pm(k)-\big(\rho_n^\pm(k)\big)^{-1}} \left\{
\begin{aligned}
&\vp_{n,1}^\pm(x,k)\vp_{n,2}^\pm(t,k),&&t>x,
\\
&\vp_{n,1}^\pm(t,k)\vp_{n,2}^\pm(x,k),&&t<x.
\end{aligned} \right.
\end{gather*}

\begin{lemma}\label{lm1.2}
Let $k\in \mathbb{C}$ be small enough. Then
\begin{enumerate}\def\theenumi{(\arabic{enumi})}
\item\label{it1lm1.2}
For any function $f\in L_2(\mathbb{R};Q)$ and $k\not=0$ the
solution to the equation
\begin{equation}\label{6.3b}
\left(-\frac{d}{dx}p\frac{d}{dx}+q-\mu_n^\pm\pm
k^2\right)u=f,\quad x\in \mathbb{R},
\end{equation}
satisfying the conditions (\ref{7.10}) is given by
$u(x,k)=(\mathcal{G}_n^\pm(k)f)(x)\in\Hloc^2(\mathbb{R})$. For
$x\not\in Q$ the function $u$ is of the form
\begin{equation}\label{1.15a}
u(x,k)=\frac{\tau_n^\pm\E^{-\k_n^\pm(k)x}\Phi_{n,2}^\pm(x,k)
}{\rho_n^\pm(k)-\big(\rho_n^\pm(k)\big)^{-1}} \int\limits_Q
\vp_{n,1}^\pm(t,k)f(t)\di t,
\end{equation}
if $x$ lies to right w.r.t. $Q$, and
\begin{equation}\label{1.15b}
u(x,k)=\frac{\tau_n^\pm\E^{\k_n^\pm(k)x}\Phi_{n,1}^\pm(x,k)}{
\rho_n^\pm(k)-\big(\rho_n^\pm(k)\big)^{-1}} \int\limits_Q
\vp_{n,2}^\pm(t,k)f(t)\di t,
\end{equation}
if $x$ lies to the left w.r.t. $Q$.

\item\label{it2lm1.2}
For any $\a_1,\a_2\in \mathbb{R}$ the operator
$\mathcal{G}_n^\pm(k): L_2(Q)\to\H^2(\a_1,\a_2)$ is boundedly
meromorphic w.r.t. $k$. The representation
\begin{equation*}
\mathcal{G}_n^\pm(k)=k^{-1}\mathcal{G}_{n,-1}^\pm+
\mathcal{G}_{n,0}^\pm+k\mathcal{G}_{n,1}(k),
\end{equation*}
holds true, where the operator $\mathcal{G}_{n,-1}^\pm$ is
defined as
\begin{equation*}
(\mathcal{G}_{n,-1}^\pm f)(x)=\pm
\frac{(f,\phi_{n}^\pm)_{L_2(Q)}}{2\sqrt{|\overset{\,\textbf{.}}{D}(\mu_n^\pm)|}}
\phi_{n}^\pm(x),
\end{equation*}
while the operator $\mathcal{G}_{n,1}^\pm: L_2(Q)\to
\H^2(\a_1,\a_2)$ is bounded and holomorphic in $k$.
\end{enumerate}
\end{lemma}

\begin{remark}\label{rm7.1}
We remind that the functions $\Phi_{n,i}^\pm$ in (\ref{1.15a}),
(\ref{1.15b}) in Item~\ref{it1lm1.2} of the lemma were defined
by (\ref{1.12}), the operator $\mathcal{G}_{n,0}^\pm$ was
defined by (\ref{2.8a}), and $\phi_n^\pm$ are the eigenfunctions
of (\ref{1.3}) normalized by (\ref{1.3c}).
\end{remark}

Item~\ref{it1lm1.2} of the lemma follows directly from the
definition of the function $G_n^\pm$ and the identities
(\ref{1.13a}), (\ref{1.12}). The validity of Item~\ref{it2lm1.2}
is due to the definition of  $G_n^\pm$, holomorphy in $k$ of
$\vp_{n,i}^\pm$, $\Phi_{n,i}^\pm$ and $\rho_n^\pm$, and the
identity (\ref{1.9}).

The uniform boundedness of the operator $\mathcal{L}_\e$ and
Item~\ref{it2lm1.2} of Lemma~\ref{lm1.2} follows that for all
$\e$ and $k\in \mathbb{C}$ small enough the operator
$\mathcal{L}_\e\big(\mathcal{G}_{n,0}^\pm+
k\mathcal{G}_{n,1}^\pm(k)\big): L_2(Q)\to L_2(Q)$ is bounded
uniformly in $\e$ and $k$. Thus, for all $\e$ and $k$ small
enough the bounded operator
\begin{equation*}
\mathcal{A}_n^\pm(\e,k):=\Big(\I-\e
\mathcal{L}_\e\big(\mathcal{G}_{n,0}^\pm+
k\mathcal{G}_{n,1}^\pm(k)\big)\Big)^{-1}
\end{equation*}
is well-defined in $L_2(Q)$. Item~\ref{it2lm1.2} of
Lemma~\ref{lm1.2} follows that for all $\e$ and $k$ small enough
the operator $\mathcal{A}_n^\pm(\e,k)$ is boundedly holomorphic
w.r.t. $k$, and an uniform in $k$ convergence
\begin{equation}\label{2.3}
\mathcal{A}_n^\pm(\e,k)\xrightarrow[\e\to0]{}\I
\end{equation}
holds true.

\begin{lemma}\label{lm2.2}
For all $\e$ and $k$ small enough the equation
\begin{equation}\label{2.4}
k\mp\frac{\e}{2\sqrt{|\overset{\,\textbf{.}}{D}(\mu_n^\pm)|}}\big(
\mathcal{A}_n^\pm(\e,k)\mathcal{L}_\e\phi_{n}^\pm,
\phi_{n}^\pm\big)_{L_2(Q)}=0
\end{equation}
has the unique solution $k_{\e,n}^\pm$. The asymptotic formulas
\begin{align}
&k_{\e,n}^\pm=\pm\frac{\e}{2\sqrt{|\overset{\,\textbf{.}}{D}(\mu_n^\pm)|}}
\big( \mathcal{A}_n^\pm(\e,0)\mathcal{L}_\e\phi_{n}^\pm,
\phi_{n}^\pm \big)_{L_2(Q)} (1+\Odr(\e^2)),\label{2.5}
\\
&k_{\e,n}^\pm=\e\big(k_{n,\e}^{\pm,1}+\e
k_{n,\e}^{\pm,2}\big)+\Odr(\e^3),\label{2.6}
\end{align}
holds true, where $k_{n,\e}^{\pm,i}$ is from (\ref{2.8b}).
\end{lemma}

\begin{proof}
Since the operator $\mathcal{A}_n^\pm(\e,k)$ is holomorphic in
$k$, the function $k\mapsto \big(\phi_{n}^\pm,
\mathcal{A}_n^\pm(\e,k)\mathcal{L}_\e\phi_{n}^\pm\big)_{L_2(Q)}$
is holomorphic in $k$ for each value of $\e$. Moreover, by
(\ref{2.3}) and the uniform boundedness of $\mathcal{L}_\e$ this
function is bounded uniformly in $\e$ and $k$. Let $\d$ be a
small number. Then for $\e$ small enough and $|k|=\d$ we have
the estimate
\begin{equation*}
\e\left|\frac{\pm
1}{2\sqrt{|\overset{\,\textbf{.}}{D}(\mu_n^\pm)|}}\big(
\mathcal{A}_n^\pm(\e,k)\mathcal{L}_\e\phi_{n}^\pm, \phi_{n}^\pm
\big)_{L_2(Q)} \right|<|k|.
\end{equation*}
By Rouche theorem it follows that the function
\begin{equation*}
k\mapsto
k\mp\frac{\e}{2\sqrt{|\overset{\,\textbf{.}}{D}(\mu_n^\pm)|}}\big(
\mathcal{A}_n^\pm(\e,k)\mathcal{L}_\e\phi_{n}^\pm, \phi_{n}^\pm
\big)_{L_2(Q)}
\end{equation*}
has the same amount of zeros inside the disk $|k|\leqslant\d$ as
the function $k\mapsto k$ does, i.e., exactly one zero. Thus,
for $\e$ small enough the equation (\ref{2.4}) has the unique
root inside the disk $|k|\leqslant\d$ which we indicate as
$k_{\e,n}^\pm$.

An obvious identity
\begin{equation*}
\frac{\p \mathcal{A}_n^\pm}{\p k}\Big(\I-\e
\mathcal{L}_\e\big(\mathcal{G}_{n,0}^\pm+
k\mathcal{G}_{n,1}^\pm\big)\Big)-\e \mathcal{A}_n^\pm
\mathcal{L}_\e\left(\mathcal{G}_{n,0}^\pm+
k\frac{\mathcal{G}_{n,1}^\pm}{\p k}\right)=0,
\end{equation*}
implies
\begin{equation*}
\frac{\p \mathcal{A}_n^\pm}{\p k}=\e \mathcal{A}_n^\pm
\mathcal{L}_\e\left(\mathcal{G}_{n,0}^\pm+
k\frac{\mathcal{G}_{n,1}^\pm}{\p k}\right)\mathcal{A}_n^\pm.
\end{equation*}
This formula, the convergence (\ref{2.3}) and
Item~\ref{it2lm1.2} of Lemma~\ref{lm1.2} give rise to the
uniform in $\e$ and $k$ estimate
\begin{equation*}
\Big\|\frac{\p \mathcal{A}_n^\pm}{\p k}\Big\|\leqslant C\e.
\end{equation*}
Employing this estimate and the formula
\begin{equation*}
\mathcal{A}_n^\pm(\e,k)-\mathcal{A}_n^\pm(\e,0)=\int\limits_0^k
\frac{\p \mathcal{A}_n^\pm}{\p k}(\e,z)\di z,
\end{equation*}
we obtain that
\begin{equation*}
\mathcal{A}_n^\pm(\e,k)=\mathcal{A}_n^\pm(\e,0)+\Odr(\e |k|).
\end{equation*}
We substitute this identity into (\ref{2.4}) to get
\begin{gather*}
k_{\e,n}^\pm\mp\frac{\e}{2\sqrt{|\overset{\,\textbf{.}}{D}(\mu_n^\pm)|}}\big(
\mathcal{A}_n^\pm(\e,0)\mathcal{L}_\e\phi_{n}^\pm, \phi_{n}^\pm
\big)_{L_2(Q)} +\Odr(\e^2|k_{\e,n}^\pm|)=0,
\\
k_{\e,n}^\pm\big(1+\Odr(\e^2)\big)=\pm\frac{\e}{2\sqrt{|\overset{\,\textbf{.}}{D}(\mu_n^\pm)|}}
\big( \mathcal{A}_n^\pm(\e,0)\mathcal{L}_\e\phi_{n}^\pm,
\phi_{n}^\pm \big)_{L_2(Q)}.
\end{gather*}
The last relation implies the asymptotics (\ref{2.5}).

Due to the definition of $\mathcal{A}_n^\pm(\e,k)$ the identity
\begin{equation}\label{3.1}
\mathcal{A}_n^\pm(\e,k)=\I+\e \mathcal{L}_\e
\mathcal{G}_{n,0}^\pm+\Odr(\e^2+\e |k|)
\end{equation}
holds true. We substitute this identity into the asymptotic
(\ref{2.5}) and deduce that
\begin{equation*}
k_{\e,n}^\pm=\e \big(k_{n,\e}^{\pm,1}+\e
k_{n,\e}^{\pm,2}\big)+\Odr\big(\e^3+\e^2 |k_{\e,n}^\pm|\big).
\end{equation*}
Now it sufficient to employ the estimate $k_{\e,n}^\pm=\Odr(\e)$
yielded by (\ref{2.5}) to complete the proof.
\end{proof}

\section{Proof of Theorems~\ref{th2.5},~\ref{th2.6} }

\begin{proof}[Proof of Theorem~\ref{th2.5}]
In the proof we employ the approach suggested in \cite{G1} which
is a modification of Birman-Schwinger principle. Let
$\l=\mu_n^\pm\mp k^2$ be an eigenvalue of $\mathcal{H}_\e$ lying
in a vicinity of $\mu_n^\pm$. Then an associated eigenfunction
$\psi$ satisfies the equation (\ref{6.3b}) with $f=\e
\mathcal{L}_\e\psi$ as well as to the conditions (\ref{7.10}),
where $\RE \k_n^\pm(k)>0$. By (\ref{1.10}) the last inequality
is equivalent to $\RE k>0$. By Item~\ref{it1lm1.2} of
Lemma~\ref{lm1.2} we conclude that $\psi=\mathcal{G}_n^\pm(k)f$.
We apply the operator $\e \mathcal{L}_\e$ to this relation to
obtain the equation for the function $f$,
\begin{equation}\label{6.2}
f-\e \mathcal{L}_\e\mathcal{G}_n^\pm(k)f=0.
\end{equation}
The operator $\mathcal{L}_\e\mathcal{G}_n^\pm(k)$, as it follows
from Item~\ref{it1lm1.2} of Lemma~\ref{lm1.2} and uniform
boundedness of the operator $\mathcal{L}_\e$, is a bounded
operator in $L_2(Q)$. The equation (\ref{6.2}) can be hence
considered as an equation in this space.

If $f$ is a nontrivial solution to (\ref{6.2}) for some $k$, it
follows that the function $\psi:=\mathcal{G}_n^\pm(k)f$ is a
nontrivial solution to
\begin{equation}\label{6.1}
\left(-\frac{d}{dx}p\frac{d}{dx}+q-\e \mathcal{L}_\e-
\mu_n^\pm\pm k^2\right)\psi=0,\quad x\in\mathbb{R},
\end{equation}
behaving at infinity in accordance with (\ref{1.15a}),
(\ref{1.15b}). This function is an element of
$\H^2(\mathbb{R})$, if and only if $\RE k>0$. Thus, the number
$\l=\mu_n^\pm\mp k^2$ is an eigenvalue of $\mathcal{H}_\e$, if
and only if $\RE k>0$. Therefore, the problem on eigenvalues of
$\mathcal{H}_\e$ tending to $\mu_n^\pm$ as $\e\to0$ reduces to
the problem on finding the values of $k$ with $\RE k>0$ for
which the equation (\ref{6.2}) has a nontrivial solution.

In accordance with Item~\ref{it2lm1.2} of Lemma~\ref{lm1.2}, the
equation~(\ref{6.2}) can be rewritten as
\begin{equation*}
f-\frac{\e}{k}\mathcal{L}_\e\mathcal{G}_{n,-1}^\pm f-
\e\mathcal{L}_\e\big(\mathcal{G}_{n,0}^\pm+k\mathcal{G}_{n,1}(k)\big)
f=0.
\end{equation*}
We apply the operator $\mathcal{A}_n^\pm(\e,k)$ to this equation
and obtain
\begin{equation}\label{6.3}
f\mp\frac{\e(f,\phi_{n}^\pm)_{L_2(Q)}}
{2k\sqrt{|\overset{\,\textbf{.}}{D}(\mu_n^\pm)|}}
\mathcal{A}_n^\pm(\e,k)\mathcal{L}_\e\phi_{n}^\pm=0.
\end{equation}
Let $f$ be a nontrivial solution to (\ref{6.2}). In this case
$(f,\phi_{n}^\pm)_{L_2(Q)}\not=0$, since otherwise by
(\ref{6.3}) we have $f\equiv 0$. Taking this fact into account,
we calculate the inner product of equation (\ref{6.3}) with
$\phi_{n}^\pm$ in $L_2(Q)$ and arrive at the equation. Hence,
the values of the parameter $k$ for those the equation
(\ref{6.1}) has a nontrivial solution are the roots of the
equation (\ref{2.4}). For $k=k_{\e,n}^\pm$ the equation
(\ref{6.3}), and, therefore, (\ref{6.2}) has a solution
\begin{equation*}
f_{\e,n}^\pm=\e
\mathcal{A}_n^\pm(\e,k_{\e,n}^\pm)\mathcal{L}_\e\phi_{n}^\pm.
\end{equation*}
This fact can be checked easily by substituting $f_{\e,n}^\pm$
into (\ref{6.3}) and bearing in mind (\ref{2.4}). The
corresponding solution of the equation (\ref{6.1}) is given by
\begin{equation}\label{2.12}
\psi_{\e,n}^\pm=\e \mathcal{G}_n^\pm\big(k_{\e,n}^\pm\big)
\mathcal{A}_n^\pm\big(\e,k_{\e,n}^\pm\big)
\mathcal{L}_\e\phi_{n}^\pm.
\end{equation}
This solution is nontrivial since by Item~\ref{it2lm1.2} of
Lemma~\ref{lm1.2}, (\ref{3.1}), (\ref{2.4}) and uniform in $\e$
boundedness of $\mathcal{L}_\e$ for any $\a_1,\a_2\in
\mathbb{R}$ in the norm of $W_2^2(\a_1,\a_2)$ the identity
(\ref{2.13}) is valid. Therefore, $f_{\e,n}^\pm\not\equiv0$. The
relation (\ref{6.3}) implies that all the solutions to
(\ref{6.2}) with $k=k_{\e,n}^\pm$ are proportional to
$f_{\e,n}^\pm$, and this is why the solution $f_{\e,n}^\pm$ is
unique up to a multiplicative constant. Thus for $\e$ small
enough there exists the unique $k$ in a small neighbourhood of
zero for which the equation (\ref{6.2}) has a nontrivial
solution. This value of $k$ is a root of the equation
(\ref{2.4}). The function $\psi_{\e,n}^\pm$ is an eigenfunction
of the operator $\mathcal{H}_\e$ only in the case $\RE
k_{\e,n}^\pm>0$. Otherwise the operator  $\mathcal{H}_\e$ has no
eigenvalues converging to $\mu_n^\pm$ as $\e\to0$. Therefore,
the operator $\mathcal{H}_\e$ has at most one eigenvalue which
converges to $\mu_n^\pm$. If exists, it is simple and is given
by $\l_{\e,n}^\pm=\mu_n^\pm\mp \big(k_{\e,n}^\pm\big)^2$. By
(\ref{2.5}) we obtain that the inequality $\RE k_{\e,n}^\pm>0$
is equivalent to (\ref{2.9}). If this inequality holds true, the
asymptotic expansions for $\l_{\e,n}^\pm$ follows immediately
from (\ref{2.5}), (\ref{2.6}). The associated eigenfunction
satisfies the expansion (\ref{2.13}).
\end{proof}

\begin{proof}[Proof of Theorem~\ref{th2.6}]
As it was established in the proof of Theorem~\ref{th2.5}, the
criterion of the existence of the eigenvalue is the inequality
$\RE k_{\e,n}^\pm>0$. The asymptotic (\ref{2.6}) implies that
the sufficient condition this inequality to be true is the
estimate (\ref{1.12a}), while the sufficient condition of
violation is the estimate (\ref{1.12b}).
\end{proof}

In conclusion let us find out the behaviour of $\psi_{\e,n}^\pm$
at infinity. This function satisfies (\ref{1.15a}),
(\ref{1.15b}), and hence for $x$ lying to the right w.r.t. $Q$
we have the relation
\begin{align*}
&\psi_{\e,n}^\pm(x)=c(\e)\E^{-\k_n^\pm(k_{\e,n}^\pm)}\Phi_{n,2}^\pm
(x,k_{\e,n}^\pm),
\\
&c(\e)=\frac{\tau_n^\pm}{\rho_n^\pm(k_{\e,n}^\pm)-\big(
\rho_n^\pm(k_{\e,n}^\pm)\big)^{-1}}\int\limits_\mathbb{R}
\vp_{n,1}^\pm(t,k_{\e,n}^\pm)g_{\e,n}^\pm(t)\di t.
\end{align*}
The identity
$\vp_{n,1}^\pm(\cdot,k_{\e,n}^\pm)=\phi_{n}^\pm(\cdot)+
\Odr(|k_{\e,n}^\pm|)$, the asymptotics (\ref{1.9}) and the
equation (\ref{2.4}) yield
\begin{align*}
c(\e)&=\pm\frac{\e}{2\sqrt{|\overset{\,\textbf{.}}{D}(\mu_n^\pm)|}k_{\e,n}^\pm\big(1+
\Odr\big(|k_{\e,n}|\big)\big)}\left(
\big(\phi_{n}^\pm,\mathcal{A}_n^\pm(\e,k_{\e,n}^\pm)\mathcal{L}_\e
\phi_{n}^\pm\big)_{L_2(Q)}+\Odr\big(|k_{\e,n}^\pm|\big) \right)=
\\
&=\left(\pm\frac{\e}{2\sqrt{|\overset{\,\textbf{.}}{D}(\mu_n^\pm)|}k_{\e,n}^\pm}
\big(\phi_{n}^\pm,\mathcal{A}_n^\pm(\e,k_{\e,n}^\pm)\mathcal{L}_\e
\phi_{n}^\pm\big)_{L_2(Q)}+\Odr\big(\e
|k_{\e,n}^\pm|\big)\right)
\left(1+\Odr\big(|k_{\e,n}^\pm|\big)\right)
=
\\
&=1+\Odr\big(|k_{\e,n}^\pm|\big).
\end{align*}
Thus,
\begin{equation}\label{2.14}
\psi_{\e,n}^\pm(x)=\big(1+\Odr(
|k_{n,\e}^\pm|)\big)\E^{-\k_n^\pm(k_{n,\e}^\pm)x}
\Phi_{n,2}^\pm(x,k_{n,\e}^\pm),
\end{equation}
if $x$ lies to the right w.r.t. $Q$. By analogy one can prove
that
\begin{equation}\label{2.15}
\psi_{\e,n}^\pm(x)=\big(1+\Odr(
|k_{n,\e}^\pm|)\big)\E^{\k_n^\pm(k_{n,\e}^\pm)x}
\Phi_{n,1}^\pm(x,k_{n,\e}^\pm),
\end{equation}
if $x$ lies to the left w.r.t. $Q$.

\section{Examples}

In this section we give some examples of the operator
$\mathcal{L}_\e$. Throughout the section the symbol $Q$
indicates certain fixed finite interval.

\textbf{1. Second order differential operator.} Let
\begin{equation*}
\mathcal{L}_\e:=b_2(x,\e)\frac{d^2}{dx^2}+b_1(x,\e)\frac{d}{dx}
+b_0(x,\e),
\end{equation*}
where $b_i\in L_\infty(Q)$ are complex-valued functions so that
$\supp b_i(\cdot,\e)\subseteq Q$, and the norms
$\|b_i(\cdot,\e)\|_{L_\infty(Q)}$ are bounded uniformly in $\e$.
The operator $\mathcal{L}_\e$ satisfies the estimate
(\ref{1.3e}), and hence by Theorem~\ref{th2.31} the continuous
spectrum of $\mathcal{H}_\e$ contains no embedded eigenvalues.
The primitives of the functions $b_i$, $q$, $p$, and the
function $p$ are continuous and have bounded variation (see, for
instance, \cite[Ch. V\!I, Sec. 1,2]{KF}). Employing this fact
and applying the uniqueness theorem for Cauchy problem from
\cite[Ch. 1, Sec. 3, Item 3]{F}, one can check easily that each
eigenvalue of $\mathcal{H}_\e$ is simple. The existence and the
asymptotics of the eigenvalues converging to the edges of the
non-degenerate lacunas are described by
Theorems~\ref{th2.5},~\ref{th2.6}.

In the particular case $b_1=b_2\equiv0$ the coefficient
$k_{\e,n}^{\pm,1}$ reads as follows
\begin{equation}\label{4.2}
k_{n,\e}^{\pm,1}=\pm\frac{(\phi_n^\pm,
b_0\phi_n^\pm)_{L_2(Q)}}{2\sqrt{|\overset{\,\textbf{.}}{D}(\mu_n^\pm)|}}.
\end{equation}
If $b_0$ is independent of $\e$ and $b_0\not\equiv0$,
$b_0\geqslant 0$, it follows that this coefficient is
independent of $\e$ and $k_{n,\e}^{+,1}>0$, $k_{n,\e}^{-,1}<0$.
In this case we can employ Theorem~\ref{th2.6}, where the
estimate (\ref{1.12a}) is valid for $k_{n,\e}^{+,1}+\e
k_{n,\e}^{+,2}$ with $C(\e)=\e^{-1}$, and the estimate
(\ref{1.12b}) does for $k_{n,\e}^{-,1}+\e k_{n,\e}^{-,2}$ with
$C(\e)=\e^{-1}$. Therefore, the eigenvalues exist only near the
right edges of non-degenerate lacunas and their asymptotics are
due to (\ref{2.11}). If $b_0\not\equiv0$ and $b_0\leqslant 0$,
we have $k_{n,\e}^{+,1}<0$, $k_{n,\e}^{-,1}>0$, and in this case
the eigenvalues exist only near the right edges of the
non-degenerate lacunas. By (\ref{2.11}) and (\ref{4.2}) in this
particular case the asymptotics of these eigenvalues are as
follows,
\begin{equation*}
\l_{\e,n}^\pm=\mu_n^\pm\mp \e^2 \frac{(\phi_n^\pm,
b_0\phi_n^\pm)_{L_2(Q)}^2}{4|\overset{\,\textbf{.}}{D}(\mu_n^\pm)|}
+\Odr(\e^3).
\end{equation*}

\textbf{2. Integral operator.} Let
\begin{equation*}
(\mathcal{L}_\e u)(x):=\int\limits_Q L_\e(x,y)u(y)dy,
\end{equation*}
where $\supp L_\e(\cdot,y)\subseteq Q$ and an uniform in $\e$
estimate
\begin{equation*}
\int\limits_{Q\times Q} |L_\e(x,y)|^2\di x\di y\leqslant C
\end{equation*}
is valid. The operator $\mathcal{L}_\e$ is  bounded uniformly in
$\e$ as an operator in  $L_2(Q)$, and  the statement of
Item~\ref{it2th2.4} of Theorem~\ref{th2.31} with $a_\e\equiv0$
thus holds true for this operator. Therefore, the continuous
spectrum of $\mathcal{H}_\e$ does not contain embedded
eigenvalues. The coefficient $k_{n,\e}^{\pm,1}$ in (\ref{2.11})
is determined by the formula
\begin{equation*}
k_{n,\e}^{\pm,1}=\pm\int\limits_{Q\times Q}
L_\e(x,y)\phi_n^\pm(x)\phi_n^\pm(y)\di x\di y.
\end{equation*}
In particular, if $L_\e(x,y)=\b \overline{b(x)}b(y)$, $\supp
b\subseteq Q$, $b\in L_2(Q)$, and $\b\in \mathbb{C}$ is a
constant, it follows that
\begin{equation}\label{9.01}
k_{n,\e}^{\pm,1}=\pm\frac{\b\left|
(b,\phi_n^\pm)_{L_2(Q)}\right|^2}{2\sqrt{|\overset{\,\textbf{.}}{D}(\mu_n^\pm)|}},
\quad k_{n,\e}^{\pm,2}=k_{n,\e}^{\pm,1} \int\limits_{Q\times Q}
G_n^\pm(x,y)\overline{b(x)}b(y)\di x\di y.
\end{equation}
In the case $(b,\phi_n^\pm)_{L_2(Q)}\not=0$ the coefficients
$k_{n,\e}^{\pm,1}$ are independent of $\e$, and
Theorem~\ref{th2.6} with $C(\e)=\e^{-1}$ is applicable. In
accordance with this theorem, the eigenvalues exist near the
right edges of the non-degenerate lacunas in the case $\RE\b>0$,
and near the left edges if $\RE\b<0$. The asymptotics expansion
for these eigenvalues are given by (\ref{2.11}) and due to
(\ref{9.01}) read as follows:
\begin{equation*}
\l_{\e,n}^\pm=\mu_n^\pm\mp \e^2 \frac{\b^2\left|
(b,\phi_n^\pm)_{L_2(Q)}\right|^4}{4|\overset{\,\textbf{.}}{D}(\mu_n^\pm)|}
\left(1+2\e\int\limits_{Q\times Q}
G_n^\pm(x,y)\overline{b(x)}b(y)\di x\di y\right)+\Odr(\e^4).
\end{equation*}

\textbf{3. Linear functional.} Let
\begin{equation*}
\mathcal{L}_\e u:=b_\e l_\e u,
\end{equation*}
where $b_\e$ is a complex-valued function such that $\supp
b_\e\subseteq Q$, and the norm $\|b_\e\|_{L_2(Q)}$ is bounded
uniformly in $\e$. The symbol $l_\e$ indicates the functional
from $\H^2(Q)$ in $\mathbb{C}$ bounded uniformly in $\e$. In
this case the operator $\mathcal{L}_\e$ defined above satisfy
all the requirements. As the example in the fifth section shows,
there exists a function $b_\e$ and a functional $l_\e$, for
which the operator $\mathcal{H}_\e$ has embedded eigenvalues. At
the same time, if $l_\e$ is a functional from $\H^1(Q)$ into
$\mathbb{C}$ bounded uniformly in $\e$, it follows that the
operator $\mathcal{L}_\e$ satisfies the hypothesis of
Item~\ref{it2th2.4} of Theorem~\ref{th2.31} with $a_\e=0$, and
in this case the operator $\mathcal{H}_\e$ has no embedded
eigenvalues.

The operator $\mathcal{L}_\e$ in this example is
finite-dimensional that allows to find the function
$\mathcal{A}_n^\pm(\e,0)\mathcal{L}_\e\phi_{n}^\pm$ explicitly:
\begin{equation*}
\mathcal{A}_n^\pm(\e,0)\mathcal{L}_\e\phi_{n}^\pm= \frac{b_\e
l_\e\phi_n^\pm}{1-\e l_\e \mathcal{G}_{n,0}^\pm b_\e}.
\end{equation*}
It follows that
\begin{equation}\label{9.02}
\big(\phi_{n}^\pm, \mathcal{A}_n^\pm(\e,0)
\mathcal{L}_\e\phi_{n}^\pm\big)_{L_2(Q)}=
\frac{(b_\e,\phi_n^\pm)_{L_2(Q)}l_\e\phi_n^\pm}{1-\e l_\e
\mathcal{G}_{n,0}^\pm b_\e}.
\end{equation}
Now by Theorem~\ref{th2.5} we obtain that the eigenvalue
$\l_{\e,n}^\pm$ exists if and only if $\pm\RE
(b_\e,\phi_n^\pm)_{L_2(Q)}l_\e\phi_n^\pm>0$. If exists, the
asymptotics for the eigenvalue $\l_{\e,n}^\pm$ is determined by
the identities (\ref{2.10}) and (\ref{9.02}):
\begin{align*}
\l_{\e,n}^\pm=&\mu_n^\pm\mp\frac{\e^2}{4|\overset{\,\textbf{.}}{D}(\mu_n^\pm)|}
\frac{\big((b_\e,\phi_n^\pm)_{L_2(Q)}l_\e\phi_n^\pm\big)^2}
{(1-\e l_\e \mathcal{G}_{n,0}^\pm
b_\e)^2}\big(1+\Odr(\e^2)\big),
\\
\l_{\e,n}^\pm=&\mu_n^\pm\mp\frac{\e^2}{4|\overset{\,\textbf{.}}{D}(\mu_n^\pm)|}
\big((b_\e,\phi_n^\pm)_{L_2(Q)}l_\e\phi_n^\pm\big)^2(1-2\e
l_\e\mathcal{G}_{n,0}^\pm b_\e)+\Odr(\e^4).
\end{align*}

\end{document}